\documentclass[12pt]{article}
\usepackage{amssymb}
\usepackage{amscd}
\usepackage{epsf}
\newtheorem{thm}{Theorem}[section]
\newtheorem{lem}[thm]{Lemma}
\newtheorem{prop}[thm]{Proposition}

\newtheorem{cor}[thm]{Corollary}
\newtheorem{defi}[thm]{Definition}

\newcommand{\inv}{^{-1}}
\newcommand{\iso}{\stackrel{\sim}{\longrightarrow}}
\newcommand{\tim}{^\times}
\newcommand{\C}{{\mathbb{C}}}
\newcommand{\Z}{{\mathbb{Z}}}
\newcommand{\R}{{\mathbb{R}}}

\newcommand{\ov}{\overline}

\newcommand{\Phired}{\Phi^{\rm red}}
\newcommand{\Phiredp}{\Phi^{+ \rm{red}}}
\newcommand{\Phiredm}{\Phi^{- \rm{red}}}
\newcommand{\Phiredpm}{\Phi^{\pm \rm{red}}}

\newcommand{\om}{\Omega}

\newcommand{\bg}{{\bf  G}}
\newcommand{\bn}{{\bf N}}
\newcommand{\bt}{{\bf T}}
\newcommand{\bz}{{\bf Z}}
\newcommand{\bs}{{\bf S}}
\newcommand{\bu}{{\bf U}}
\newcommand{\bc}{{\bf C}}
\newcommand{\bgl}{\bf GL}
\newcommand{\bsl}{\bf SL}
\newcommand{\bpgl}{\bf PGL}

\newcommand{\bl}{{\bf L}}
\newcommand{\bb}{\bf B}
\newcommand{\bp}{\bf P}
\newcommand{\bd}{\bf D}
\newcommand{\bm}{\bf M}
\newcommand{\bh}{\bf H}

\newcommand{\fyd}{F_Y^\Delta}
\newcommand{\fydp}{F_{Y'}^{\Delta'}}
\newcommand{\ftyd}{\tilde{F}_Y^\Delta}
\newcommand{\hiwei}{\lambda_0( \Delta)}

\newcommand{\kbar}{{\ov{K}}}

\newcommand{\lkbar}{\! \! ~_{\kbar}}

\newcommand{\larho}{\lambda_0^\rho(\Delta)}
\newcommand{\lasigma}{\lambda_0^\sigma(\Delta)}
\begin{document}

\title{Compactifications of Bruhat-Tits buildings associated to linear representations}
\author{Annette Werner}

\date{}

\maketitle

\centerline{\bf Abstract}
\small 
Let $\bg$ be a connected semisimple group over a non-Archimedean local field. For every faithful, geometrically irreducible linear representation of $\bg$ we define a compactification of the associated Bruhat-Tits building $X(\bg)$. This yields a finite family of compactifications of $X(\bg)$ which contains the polyhedral compactification. Besides, this family can be seen as a non-Archimedean analogue of the family of  Satake compactifications for symmetric spaces.

\centerline{{\bf 2000 MSC:} 20E42, 20G25} 

\normalsize 
\section*{Introduction}
In this paper we define compactifications of Bruhat-Tits buildings, which are analogous to Satake compactifications of symmetric spaces of non-compact type.

Let $\bg$ be a connected semisimple group over a non-Archimdedean local field $K$. Then the Bruhat-Tits building $X(\bg)$ is a complete metric space, endowed with a continuous $\bg(K)$-action and a poly-simplicial structure. It describes in a geometric way the combinatorics of the maximal compact
subgroups of $\bg(K)$. Besides, it can be seen as a non-Archimedean analogue of the symmetric spaces of non-compact type given by a semisimple Lie group.
In the case $\bg = {\bsl}(2,K)$, the building is an infinite regular tree. 
Such a tree can be compactified in a natural way by adding endpoints to all half-lines starting in a fixed vertex. 

However, in the higher rank case there are different ways of compactifying the building.
In \cite{bose}, the Borel-Serre compactification is used for the study of the cohomology of arithmetic groups. Its boundary can be identified with the Tits building of $\bg$. 
In \cite{la}, based on ideas of Gerardin  \cite{ge}, Landvogt constructs the polyhedral compactification of $X(\bg)$. Here the boundary consists of all Bruhat-Tits buildings corresponding to parabolics in $\bg$. In \cite{guire} a group-theoretic compactification of the vertex set is defined, which can be
identified with the closure of the vertex set in the polyhedral compactification.

In \cite{we1} and \cite{we2} two compactifications  in the case
$\bg = \bpgl(V)$ are investigated. These spaces have a natural description in terms of lattices respectively seminorms.
Besides, the compactification in \cite{we2} is related to Drinfeld's $p$-adic upper half space, and
it can be identified topologically with a closed subset of the projective Berkovich analytic space.

Considering these results, it is a natural question to ask how these compactifications for $ \bpgl$ are related to the polyhedral compactification and how they can be generalized to arbitrary semisimple groups. This question is answered in the present paper. We show  that  every faithful geometrically irreducible representation $\rho: \bg \rightarrow {\bgl}(V)$ on a finite-dimensional $K$-vector space $V$ 
gives rise to a compactification $\ov{X}(\bg)_\rho$ of the Bruhat-Tits building for $\bg$. Let us give an outline of this construction. Consider  a maximal $K$-split torus $\bs$  in $\bg$, and let $A$ be the apartment in $X(\bg)$ corresponding to $\bs$. In section 2, we use the combinatorics of the weights of $\rho$ to define a decomposition of $A$ into faces. This decomposition gives rise to a compactification $\ov{A}$ of the space $A$, which can furthermore be endowed with an action by $N$, the $K$-rational points of the normalizer of $\bs$. 

In section 3 we define for all  points $x$ in the compactified apartment $\ov{A}$  a group $P_x$, which will later turn out to be the stabilizer of $x$ in $\ov{X}(\bg)_\rho$. With these data, the compactified building can be defined in a similar way as the building itself, namely as the quotient of $\bg(K) \times \ov{A}$ after a certain equivalence relation involving the stabilizers $P_x$, see Definition 4.1. In order to show that the resulting space is compact, we need the mixed Bruhat decomposition Theorem 3.9 and the rather technical statement in Theorem 4.2. 
In Theorem 4.5, we show that the compactifications given by two representations $\rho$ and $\sigma$ coincide, if and only if the highest weights of $\rho$ and $\sigma$ lie in the same chamber face.
Hence our construction yields a finite familiy of building compactifications. If the highest weight of 
$\rho$ lies in the open chamber, then we retrieve Landvogt's polyhedral compactification.
Also the compactifications defined in \cite{we1} and \cite{we2} belong to our family, as we show in section 5. However, the Borel-Serre compactification is of a different nature and does not occur among our spaces.

Our results can be seen as a non-Archimedean analogue of Satake's compactifications of Riemann symmetric spaces of non-compact type. In \cite{sa}, Satake defines compactifications of such Riemann spaces ${\cal S } = G/K$ associated to certain complex representations of $G$. He fixes one concrete compactification of the symmetric space associated to $PSL(n, \C)$, and uses the representation to embed $\cal S$ in this space. Taking the closure gives a compactification of $\cal S$. Although his construction differs from ours, several properties of the resulting spaces also hold for our building compactifications. E.g., his spaces also only depend on the face containing the highest weight of the representation. In a future research project it will be shown  that the boundaries of our spaces $X(\bg)_\rho$ can be described in an analogous way as the boundaries in the Satake compactifications. Another subject of future research is a  description of the compactified buildings as subsets of Berkovich analytic spaces, which generalizes the result in \cite{we2}.

\section{Rational representations}
Let $K$ be a non-Archimedean locally compact field, and let $\mathfrak{O}_K$  the ring of integers and $\ov{K}$ an algebraic closure.

We consider a connected, semi-simple group $\bg$ over $K$ and fix a maximal $K$-split torus $\bs$ in $\bg$. Note that $\bg$ is geometrically connected by \cite{ega4}, 4.5.14.
Let $\bn$ be the normalizer and $\bz$ be the centralizer of $\bs$ in $\bg$. By $G$, $S$,  $N$, respectively $Z$ we denote the corresponding  groups of $K$-rational points. They all carry topologies inherited from $K$.

Throughout this paper we fix a faithful $K$-rational linear representation of $\bg$, i.e. an injective 
$K$-homomorphism
\[\rho: \bg \rightarrow {\bgl} (V),\]
where $V$ is a finite-dimensional $K$-vector space.
We assume that $\rho$ is geometrically irreducible, i.e. that the  base change
\[\rho_\kbar: \bg_\kbar \rightarrow {\bgl}( V_\kbar)\]
 is irreducible. Here we put $V_\kbar = V \otimes_K \kbar$ and $\bg_\kbar = G \otimes_K \kbar$.

Besides, we
fix a maximal torus $\bt$ in $\bg$ containing $\bs$. We write $\lkbar\Phi = \Phi (\bt_\kbar, G_\kbar)$ and $ \Phi = \Phi(\bs, \bg) = \Phi(\bs_\kbar, \bg_\kbar)$ for the corresponding root systems. The natural restriction map $j : X^\ast(\bt_\kbar) \rightarrow X^\ast(\bs_\kbar)$ induces a map $j: \lkbar\Phi \rightarrow  \Phi \cup \{0\}$. If $ \Delta $ is a basis for the root system $ \Phi$, then there exists a basis $\lkbar\Delta$ for $\lkbar\Phi$ such that 
 \[  \Delta \subset j(\lkbar \Delta) \subset  \Delta \cup \{0\},\]
see \cite{bo}, 21.8. 

Let $\lkbar W$ be the Weyl group of $\bt$, and let $ W = N/Z$ be the Weyl group of $\bs$.

Then the Borel group $\bb$ in $\bg_\kbar$ corresponding to the basis $\lkbar\Delta$ of $\lkbar\Phi$ stabilizes a uniquely determined line in $V_\kbar$, on which $\bt_\kbar$ acts via a character $\mu_0 = \mu_0(\lkbar \Delta) $ of $\bt_\kbar$, the highest weight.

Let $\Omega  \subset X^\ast(\bt_\kbar)$ denote the set of all weights of $\bt_\kbar$. 
For every $\mu \in \Omega$ the difference $\mu_0 - \mu$ lies in $\lkbar\Phi^+$, i.e. it is of the form
$\mu_0 - \mu = \sum_{a \in \lkbar\Delta} n_a a$ for non-negative integers $n_a$.

A $K$-weight $\lambda \in X^\ast(\bs_\kbar)$ is defined as the restriction of a weight $\mu \in \Omega$ to $\bs_\kbar$, i.e. $\lambda = j(\mu)$.  By $\lambda_0 = \lambda_0( \Delta)$  we denote the $K$-weight induced by the highest weight $\mu_0( \lkbar\Delta)$. It only depends on $ \Delta$, not on the choice of  $\lkbar\Delta$ compatible with $\Delta$. For any $K$-weight $\lambda$ we have
\[\lambda_0(\Delta) - \lambda = \sum_{a \in  \Delta} n_a a\]
with $n_a \geq 0$. We write 
$[\lambda_0 - \lambda] = \{ a \in \Delta: n_a \neq 0\} = \{ a \in \Delta: n_a > 0\}$ for the support of 
$\lambda_0 - \lambda$.

For every $K$-weight $\lambda$ we denote the corresponding weight space by $V_\lambda = \{ v \in V: \rho(s) (v) = \lambda(s) v \mbox{ for all } s \in S\}$. The vector space $V$ is the direct sum of its weight spaces $V_\lambda$. 

\begin{defi}
We associate to every  finite subset $M$ of $X^\ast(\bs)$ the graph with vertex set $M$, such that  $m$ and $n$ are connected by an edge, if and only if $(m,n) \neq 0$,  where $(,)$ denotes a $ W$-invariant scalar product on $X^\ast(\bs)$. We call $M$ connected, if the associated graph is connected. 

Let $\Delta$ be a basis of $\Phi$. A subset $Y \subset \Delta$ is called admissible, if 
the set $Y \cup \{\lambda_0(\Delta)\}$ is connected. 
\end{defi}

A similar definition in the context of real Lie groups was used in \cite{sa}.

We will frequently apply the following result 12.16 from  \cite{boti}:
A set $Y \subset \Delta$ is admissible if and only if there exists a $K$-weight $\lambda$ such that
$[\lambda_0( \Delta) - \lambda] = Y$.

\section{Compactification of one apartment}

Let $X = X(\bg)$ be the Bruhat-Tits building associated to $\bg$.
The maximal $K$-split torus $\bs$ in $\bg$ defines an apartment $A$ in $X$. In this section, we use the representation $\rho$ to define a decomposition of $A$ into faces. 

We denote by $X_\ast(\bs)$ respectively $X^\ast(\bs)$ the group of characters respectively cocharacters of $\bs$.  Then the  apartment defined by $\bs$ is  the real vector space
$A = X_\ast(\bs) \otimes \R$.
It carries a natural topology. 

We use the natural pairing $X_\ast(\bs) \times X^\ast(\bs) \rightarrow \Z$ to identify the dual space $A^\ast$ with $X^\ast(\bs) \otimes \R$. 
In particular, the root system $\Phi= \Phi(\bs, \bg)$  with Weyl group $W = N/Z$  lies in the dual space $A^\ast$ of $A$.

For every basis $\Delta$ of $\Phi$ we choose a basis $\lkbar\Delta$ of $\lkbar\Phi$ as in the preceeding section, i.e. so that $ \Delta \subset j(\lkbar \Delta) \subset  \Delta \cup \{0\}$, where $j$ is the restriction map $j : X^\ast(\bt_\kbar) \rightarrow X^\ast(\bs_\kbar)$. Then the highest weight $\mu_0( \lkbar\Delta)$ induces a $K$-weight $\lambda_0 = \lambda_0(\Delta)$.

For every admissible subset $Y \subset \Delta$ we put
\begin{eqnarray*}
F_Y^{\Delta} &= \{ x \in A: & a(x) = 0 \mbox{ for all } a \in Y \mbox{ and }\\ 
 & & (\lambda_0-\lambda)(x) > 0 \mbox{ for all $K$-weights } \lambda \mbox{ such that } [\lambda_0 - \lambda] \not\subset Y.\}
\end{eqnarray*}

All $\fyd$ are non-empty convex subsets of $A$. 

Let us denote by $\langle Y\rangle $ the real subspace of the dual space $A^\ast$ which is generated by $Y$. 

\begin{lem}
i) For bases $\Delta$ and $\Delta'$ of $\Phi$ and admissible subsets $Y \subset \Delta$ and  $Y' \subset \Delta'$ we have $\fyd = \fydp$, if and only if $\langle Y\rangle  = \langle Y'\rangle $ and $\lambda_0(\Delta) -\lambda_0(\Delta') \in \langle Y\rangle  = \langle Y'\rangle $.

i) Two sets $\fyd$ and $\fydp$ are either equal or disjoint. 

ii) For every $x\in A$ there is an admissible subset $Y$ of some $\Delta$ such that $x \in \fyd$. 
\end{lem}

{\bf Proof: }Assume that $\fyd \cap \fydp$ contains a point $x$.  Then for all $K$-weights $\lambda$ we have $(\lambda_0(\Delta) - \lambda)(x) \geq 0$, and $(\lambda_0(\Delta')- \lambda)(x) \geq 0$, hence
$(\lambda_0(\Delta') - \lambda_0(\Delta)) (x) = 0$, which implies that $\lambda_0(\Delta') - \lambda_0(\Delta)$ lies in $\langle Y'\rangle $. 

Now we show that $Y$ is contained in $\langle Y'\rangle $. We write $Y = \{ a_1, \ldots a_s\}$ in such a way that all subsets $\{a_1, \ldots, a_i\}$ are also admissible. 
By \cite{boti}, 12.16, there are $K$-weights $\lambda_i$ such that $[\lambda_0(\Delta) - \lambda_i] = \{a_1, \ldots, a_i\}.$
Assume that $a_1, \ldots, a_{k-1}$  are contained in $\langle Y' \rangle$. (Here $k = 1$ is permitted, then the condition is empty.) We can write $\lambda_0(\Delta) - \lambda_{k} = b + n a_k$ for some element $b \in \langle Y'\rangle $ and $n>0$. 
Since $(\lambda_0(\Delta) - \lambda_k)(x) = 0$, we find $(\lambda(\Delta') - \lambda_k)(x) = 0$, 
hence $[\lambda(\Delta') - \lambda_k] \subset Y'$, which implies $\lambda_0(\Delta) - \lambda_k =
(\lambda_0 (\Delta)- \lambda_0(\Delta')) + (\lambda_0(\Delta') - \lambda_k) \in \langle Y'\rangle $. Therefore
$a_k \in \langle Y'\rangle $. By induction we conclude $Y \subset \langle Y'\rangle $. 
Exchanging the roles of $Y$ and $Y'$ we find $Y \subset \langle Y'\rangle $, hence $\langle Y\rangle  = \langle Y'\rangle $.

Assume on the other hand, that $\langle Y\rangle  = \langle Y'\rangle $ and that $\lambda_0( \Delta) -  \lambda_0(\Delta') \in \langle Y \rangle$. 
Let $x$ be a point in $\fyd$. Then $b(x) = 0$ for all $b \in \langle Y\rangle $, in particular for the $b \in Y'$. Besides, if $[\lambda_0(\Delta') - \lambda] \not\subset Y'$, then $\lambda_0(\Delta') - \lambda$ is not contained in $\langle Y'\rangle  = \langle Y \rangle$. Since $\lambda_0(\Delta) -\lambda_0(\Delta') \in \langle Y \rangle$, it follows that  $[\lambda_0( \Delta) - \lambda] \not\subset Y$, which implies
that $(\lambda_0(\Delta')- \lambda)(x) = (\lambda(\Delta) - \lambda)(x) > 0$. Therefore we find $x \in \fydp$, i.e. we have $\fyd \subset \fydp$. The same argument with reversed roles shows $\fydp \subset \fyd$, so that $\fyd = \fydp$. 
This proves i) and ii).

iii) Let $x$ be a point in $A$. Then there is a basis $\Delta$ such that $x$ lies in the closed chamber
$\{ y \in A: a(y) \geq 0 \mbox{ for all }a \in \Delta\}$. Let $Y$ be the maximal admissible subset of the
set $\{a \in \Delta: a(x) = 0\}$. Since for all $K$-weights $\lambda$ the set $[\lambda_0(\Delta) - \lambda]$ is admissible, it is either contained in $Y$ or it contains a root $a \in \Delta$ such that $a(x) > 0$. Hence
$x$ lies in $\fyd$.\hfill$\Box$

Hence the sets $\fyd$ provide a decomposition of $A$ in convex subsets (faces). Let 
\[\Sigma = \{ \fyd: \Delta \mbox{ basis of }\Phi, Y \subset \Delta \mbox{ admissible}\}\]
be the collection of all these faces. 

If $\lambda_0(\Delta)$ lies in the interior of the chamber in $X^\ast(\bs)_\R$ defined by $\Delta$,
then all subsets of $\Delta$ are admissible. In this case, $\Sigma$ coincides with the set
of chamber faces in $A$ defined by the root system $\Phi$. Note however that the decomposition of $A$ into the faces  in $\Sigma$ is in general not induced by a set of hyperplanes (as investigated in \cite{bou}, chapter V). This can be seen e.g. in the case of the identical representation for $SL(3)$, cf. \cite{we1}.

By $\ftyd$ we denote the closure of $\fyd$ in A. Note that
\[\ftyd = \{x \in A: a(x) = 0 \mbox{ for all } a \in Y \mbox{ and } (\lambda_0(\Delta) - \lambda)(x) \geq 0 \mbox{ for all $K$-weights }\lambda\}. \]
In fact, consider a point y in the set on the right hand side and some $x \in \fyd$. Then every point on the line segment beween $x$ and $y$ apart from $y$ lies in $\fyd$, so that $y \in \ftyd$. The other inclusion is clear. 

The same argument as in the proof of the preceeding lemma shows that 
\[ \fydp \subset \ftyd \mbox{ iff }  [\lambda_0(\Delta') - \lambda_0(\Delta)] \subset Y' \mbox { and } \langle Y \rangle \subset \langle Y'\rangle .\]

Note that the subspace of $A$ spanned by $\fyd$ is equal to 
\[ \langle \fyd \rangle  = \{ x \in A: a(x) = 0 \mbox{ for all }a \in Y\}.\]

Now we define the boundary components of $A$.
\begin{defi}
For every face $F \in \Sigma$ let $A_F$ be the $\R$-vector space $A_F = A / \langle F \rangle$,
and let $r_F : A \rightarrow A_F$ be the quotient map.  

Put $\ov{A} = \bigcup_{F \in \Sigma} A_F$
\end{defi}

Note that the apartment $A$ is contained in $\ov{A}$. To see this, it suffices to check that $\Delta$ is admissible.
Since $F^\Delta_\Delta = \{0\}$, the corresponding quotient space is $A$ itself. By \cite{boti}, 4.28 and 5.10, the root system $\Phi$ is the direct sum of root systems $\Phi_i$ corresponding to quasi-simple subgroups $\bg_i$ of positive $K$-rank. 
Hence  $\Delta$ is the union of bases $\Delta_i$ in $\Phi_i$. If $\lambda_0(\Delta)$ is perpendicular to $\Delta_i$, then $\rho$ maps the maximal $K$-split torus in $\bg_i$ to the center of ${\bgl}(V)$. Since $\rho$ is injective, it follows that the maximal $K$-split torus in $\bg_i$ is central, which contradicts the fact that $\bg_i$ has positive $K$-rank. 
Therefore $\lambda_0(\Delta)$ is not perpendicular to $\Delta_i$. Since $\Delta_i$ is connected by \cite{boti}, 5.11 and \cite{bou}, VI, 1.7, Corollary 5, we find that  $\Delta_i$ is admissible. Hence $\Delta$ is indeed admissible. 

For all bounded open subsets $U \subset A$ and all faces $F \in \Sigma$ we define
\[\Gamma_F(U) = \bigcup_{F' \subset \tilde{F}} r_{F'} (U + F).\]
The subset $\Gamma_F(U)$  of $\ov{A}$ meets by definition exactly the components $A_{F'}$ with $F' \subset \tilde{F}$. For $F = \fyd$ these are the components $A_{F'}$ for all  $F' = \fydp$ such that $[\lambda_0(\Delta') - \lambda_0(\Delta) ] \subset Y'$ and $\langle Y \rangle \subset \langle Y' \rangle$. In particular, every $\Gamma_F(U)$ meets $A$, and we have $\Gamma_F(U) \cap A = U + F$. 

\begin{lem}
i) If $U_1 + F_1 $ and $U_2 + F_2$ are disjoint, then $\Gamma_{F_1}(U_1)$ and $\Gamma_{F_2}(U_2)$ are disjoint.

ii) The sets $\Gamma_F(U)$ form a basis of a topology on $A$.
\end{lem}

{\bf Proof: } i) Write $F_1 = F^{\Delta_1}_{Y_1}$ and $F_2 = F^{\Delta_2}_{Y_2}$, and  assume that for some face $F= \fyd$ with $F \subset \tilde{F}_1$ and $F \subset \tilde{F}_2$ there are points $x \in U_1 +F_1$ and
$y \in U_2+F_2$ such that $r_F(x) = r_F(y)$, i.e. such that $x-y \in \langle F \rangle$. 
Now $F$ contains points $f$ such that for all $\lambda$ with $[\lambda_0(\Delta) - \lambda]\not\subset Y$ the value $(\lambda_0(\Delta)-\lambda)(f)$ is arbitrarily large. Hence there is some $f \in F$ such that for all those $\lambda$ we have $(\lambda_0(\Delta) - \lambda)(x-y+f) >0$. 
Hence $x-y+f$ lies in $F$, so that $x+f = y+ (x-y+f) \in U_2 + F_2 + F \subset U_2 + \tilde{F}_2 = U_2 + F_2$. Besides,  $x + f \in U_1  + F_1 + F \subset U_1 + \tilde{F}_1 =  U_1 + F_1$, so that $U_1+F_1$ and $u_2 + F_2$ are not disjoint.

ii) Let $z$ be a point in the intersection of $\Gamma_{F_1}(U_1)$ and $\Gamma_{F_2}(U_2)$. We have to show that there exists some set $\Gamma_F(U)$ containing $x$ such that $\Gamma_F(U) \subset \Gamma_{F_1}(U_1) \cap \Gamma_{F_2}(U_2)$. Let $F = \fyd$ be the face with $z \in A_F$. Then $F$ is contained in $\tilde{F}_1$ and $\tilde{F}_2$. As in i) we find some $x \in U_1+ F_1 \cap U_2+F_2$ such that $r_F(x) = z$. Let $U$ be a bounded open neighbourhood of $x$ contained in the open set $U_1+ F_1 \cap U_2+F_2$. Then $U + F$ is contained in $U_1 + F_1$ and $U_2 + F_2$, which implies $\Gamma_F(U) \subset \Gamma_{F_1}(U_1) \cap \Gamma_{F_2}(U_2)$. \hfill$\Box$

We define a topology on $\ov{A}$ by taking as a basis all sets of the form $\Gamma_F(U)$, where
$F$ is a face in $\Sigma$, and $U$ is a bounded open subset of $A$.
The induced topology on $A$ is the given one.

Note that if $F = \fyd$, 
then the dual vector space of $A_F$ can be identified with $\langle Y \rangle \subset A^\ast$.
In particular, we can apply elements of $\langle Y \rangle$ to points in $A_F$. 

Note that a sequence of points $x_n \in A$ converges to the boundary point $x \in A_F$ for $F = \fyd$ if and only if
\begin{eqnarray*}
a(x_n) \rightarrow a(x)  & \mbox{ for all }  a \in Y \mbox{ and}\\
(\lambda_0(\Delta) - \lambda)(x_n) \rightarrow \infty & \mbox{ for all $K$-weights $\lambda$ with } [\lambda_0(\Delta) - \lambda] \not\subset Y
\end{eqnarray*}

Since the Euclidean topology on $A$ is second countable and there are only finitely many faces, $\ov{A}$ is a second countable topological space, i.e. there exists a countable basis for its topology. 
\begin{prop}
The topological space $\ov{A}$ is compact.
\end{prop}
{\bf Proof: }Let us first show that $\ov{A}$ is Hausdorff. 
Assume that $x \neq y$ lie in the same
component $A_F$, where $F = \fyd$. Then there exists some $a \in Y$ such that $a(x) \neq a(y)$. 
Choose $u$ and $v$ in $A$ such that $r_F(u)=x$ and $r_F(v)=y$. Since $u$ and $v$ assume different 
values under $a$, we find bounded open neighbourhoods $U$ around $u$ and $V$ around $v$ in $A$, such that
$a(U) \cap a(V) = \emptyset$. This implies $U+F \cap V+F = \emptyset$, hence $\Gamma_F(U)$ and $\Gamma_F(V)$ are disjoint neighbourhoods of $x$ and $y$ in $\ov{A}$ by Lemma 2.3. 

Now assume that $x \in A_F$ and $y \in A_{F'}$, where $F=\fyd$ and $F' = \fydp$ such that $[\lambda_0(\Delta)-\lambda_0(\Delta')] \not\subset Y$. Then we choose a point $v$ in $A$ with $r_{F'}(v) = y$
and a bounded open neighbourhood $V$ of $v$ in $A$. Then $(\lambda_0(\Delta') - \lambda_0(\Delta))(V+F')$ is bounded below by some real number $r$. Since $[\lambda_0(\Delta)-\lambda_0(\Delta')] \not\subset Y$ there exists a point $u\in A$ with $r_F(u) = x$ and a bounded open neighbourhood $U$ of $u$ such that $\lambda_0(\Delta)- \lambda_0(\Delta') > - r$ on $U+F$. Hence $U+F$ and $V+F'$ are disjoint, so that $\Gamma_F(U)$ and $\Gamma_{F`}(V)$ are disjoint open neighbourhoods of $x$ and $y$ in $\ov{A}$ by Lemma 2.3. 

It remains to consider the case that $x \in A_F$ and $y \in A_{F'}$ with different faces $F=\fyd$ and $F' = \fydp$ such that $\lambda_0(\Delta)-\lambda_0(\Delta') \in \langle Y \rangle \cap \langle Y' \rangle$. If $\langle Y' \rangle \subset \langle Y \rangle$,  it follows from $F \neq F'$ that $\langle Y \rangle \not\subset \langle Y' \rangle$. Hence, by reversing the roles of $Y$ and $Y'$, if necessary, we can assume that $\langle Y' \rangle \not\subset \langle Y \rangle$. We claim that there exists a $K$-weight $\lambda$ such that $\lambda_0(\Delta) - \lambda \in \langle Y' \rangle$, but
$\lambda_0(\Delta) - \lambda \not\in \langle Y \rangle$. In fact, there is an admissible subset $Z$ of $Y'$ such that all but one element in $Z$ is contained in $\langle Y \rangle$. There exists a $K$-weight $\lambda$ such that $[\lambda_0(\Delta') - \lambda] = Z$.
Then we find $\lambda_0(\Delta) - \lambda =  (\lambda_0(\Delta) - \lambda_0(\Delta')) + (\lambda_0(\Delta') - \lambda) \in \langle Y' \rangle \backslash \langle Y \rangle$.
Now choose a point $v\in A$ with $r_{F'}(v) = y$ and a bounded open neighbourhood $V$ of $v$. Since $\lambda_0(\Delta) - \lambda \in \langle Y' \rangle$, it vanishes on $F'$. Hence $(\lambda_0(\Delta) - \lambda) (V+F')$ is contained in a bounded set $I=[r,s] \subset \R$. Since $[\lambda_0(\Delta) - \lambda] \not\subset Y$, there exists a point $u \in A$ with
$r_F(u) = x$ and a bounded open neighbourhood $U$ around $u$ such that $(\lambda_0(\Delta)-\lambda) (U + F) > s$. Then $U+F$ and $V+F'$ are disjoint. Applying Lemma 2.3, we deduce that $\Gamma_F(U)$ and $\Gamma_{F`}(V)$ are disjoint open neighbourhoods of $x$ and $y$ in $\ov{A}$. Hence $\ov{A}$ is indeed Hausdorff. 

Since $\ov{A}$ is second countable, it remains to show that
every sequence $(x_n)_{n \geq 1}$ in $\ov{A}$ has a convergent subsequence. 
The set of faces $\Sigma$ is finite. Therefore we can pass to a subsequence and assume that all $x_n$ are contained in the same component $A_F$ for some face $F = F^\Delta_Z$. Let $[Z]$ denote the set of all roots in $\Phi$ which can be written as a linear combination of the roots in $ Z$. Then $[Z]$ is a root system in the dual space
$\langle Z \rangle$ of $A_F$ by \cite{bou}, VI, § 1.1. After passing to a subsequence, we can assume that all $x_n$ are contained in one closed chamber in $A_F$. This chamber corresponds to a basis $Z'$ of the root system $[Z]$. Since $Z$ is also a basis
of $[Z]$, there is an element in the Weyl group of $[Z]$ mapping $Z$ to $Z'$. Since it is a product of reflections corresponding to roots in $Z$, we can lift it to an element $w$ of $ W$. Then $w(\Delta)$ is a basis of $\Phi$ containing $Z'$. Besides, we find $w(\lambda_0(\Delta)) - \lambda_0(\Delta)  \in \langle Z \rangle$. 
Since  $\langle Z \rangle = \langle Z' \rangle$, we deduce from Lemma 2.1 that  $ F^\Delta_Z =  F^{w(\Delta)}_{Z'}$. Therefore we can replace $\Delta$ by $w(\Delta)$ and $Z$ by $Z'$ and assume that all $x_n$ are contained in the closed chamber of $A_F$ defined by the basis $Z$ of $[Z]$.
Hence $a(x_n) \geq 0 $ for all $a \in Z$. 

Now we define $Y$ to be the biggest admissible subset of $\Delta$ contained in the set
$\{a \in Z: a(x_n) \mbox{ is bounded }\}$. After passing to a subsequence of $x_n$, we can assume that for all $a \in Y$ the sequence $a(x_n)$ is convergent. 
Put $F = \fyd$ and define $x \in A_{\fyd}$ as the point satisfying $a(x) = \lim_{n \rightarrow \infty} a(x_n)$ for all $a \in Y$. 
We claim that $x_n$ converges to $x$.

Let $U$ be a bounded open subset of $A$, and let $\Gamma_{F'}(U)$ be an open neighbourhood of $x$,
where $F' = \fydp$. 
Then $x = r_F( u+f)$ for some $u \in U$ and $f \in F'$. In particular $\Gamma_{F'}(U)$ meets $A_F$, which implies $F \subset \tilde{F}'$, i.e. 
$[\lambda_0(\Delta) - \lambda_0(\Delta')] \subset Y$ and $ \langle Y' \rangle \subset \langle Y \rangle$. 

For all $a \in Y'$ we have $a(x) = a(u)$, so that $a(x_n)$ converges to $a(u)$. Hence we find a sequence
of points $u_n \in U$ converging to $u$ such that $a(u_n) = a(x_n) $ for all $a \in Y'$. 
Besides, we choose points $z_n$ in $A$ with $r_F(z_n) = x_n$ such that the sequence $a(z_n)$ converges to $\infty$ for all $a \in  \Delta\backslash Z$. Note that $a(z_n) = a(x_n)$ for all $a \in Z$, in particular for all $a \in Y$, and hence also for all $a \in Y' \subset \langle Y \rangle$.  

We want to show next that $z_n - u_n \in F'=\fydp$ for $n$ big enough. Note first that 
for all $a \in Y'$ we have $a(z_n-u_n) = a(x_n) - a(u_n) = 0$. 
Now consider a $K$-weight $\lambda$ with $[\lambda_0(\Delta') - \lambda] \not\subset Y'$. 
If $\lambda_0(\Delta') - \lambda$ is contained in $\langle Y \rangle$, then it can be written as
$\lambda_0(\Delta') - \lambda = \sum_{a \in Y} r_a a$ for some $r_a \in \R$. Since for all $a \in Y$
the sequence $a(z_n)-a(u_n)$ converges to $a(x) - a(u) = a(f)$, we find that $(\lambda_0(\Delta') - \lambda)(z_n -u_n)$ converges to $(\lambda_0( \Delta')-\lambda)(f) $, which is positive, since $f \in \fydp$. 
Therefore $(\lambda_0( \Delta')-\lambda)(z_n-u_n) >0$ for $n$ big enough.
If $\lambda_0( \Delta') -\lambda  \in \langle Z \rangle$, but $\lambda_0( \Delta') -\lambda  \notin \langle Y \rangle$, it follows from $\lambda_0( \Delta') - \lambda_0( \Delta) \in \langle Y \rangle$ that the admissible subset $[\lambda_0( \Delta) - \lambda]$ of $ \Delta$ is contained in  $Z$, but not in $Y$.
 By definition of $Y$, the set $[\lambda_0(\Delta) - \lambda]$ contains some $a \in Z$ such that $a(x_n)$ is unbounded.
Since $\lambda_0( \Delta') - \lambda \in \langle Z \rangle$, we have 
\begin{eqnarray*} (\lambda_0( \Delta') - \lambda)(z_n - u_n) & = & (\lambda_0( \Delta') - \lambda) ( x_n - u_n)\\
~ & = &   (\lambda_0( \Delta') - \lambda_0( \Delta)) ( x_n - u_n) + (\lambda_0( \Delta) - \lambda)(x_n - u_n).
\end{eqnarray*}
 The first term is zero if $[\lambda_0(\Delta') - \lambda_0(\Delta)] \subset Y'$. If this is not the case, the first term is    positive for $n$ big enough by the previous step. The second term can be written as
$(\lambda_0( \Delta) - \lambda)(x_n - u_n) = \sum_{a \in Z} n_a a(x_n - u_n)$ with $n_a \geq 0$. Since
all occurring $a(u_n)$ are bounded, all $a(x_n) \geq 0$ with at least one of them unbounded where $n_a >0$, we conclude that this
sequence goes to $\infty$. Hence for $n$ big enough we find $(\lambda_0( \Delta') - \lambda)(z_n - u_n) > 0$. 

Now consider the case that $\lambda_0( \Delta') - \lambda \notin \langle Z \rangle$. Then
$(\lambda_0( \Delta') - \lambda)(z_n - u_n) = (\lambda_0( \Delta') - \lambda_0( \Delta)) ( z_n - u_n) + (\lambda_0( \Delta) - \lambda) ( z_n - u_n)$. Again, the first term is non-negative for $n$ big enough. 
The second term can be written as $(\lambda_0( \Delta) - \lambda)(z_n - u_n) = \sum_{a \in  \Delta} n_a a(z_n - u_n)$ with $n_a \geq 0$. 
All $a(u_n)$ are bounded, and for all $a \in Z$ we have $a(z_n) = a(x_n) \geq 0$. Since $\lambda_0( \Delta') - \lambda \notin \langle Z \rangle$, we find $[\lambda_0( \Delta)- \lambda] \not\subset Z$. Hence there exists one $a \notin Z$ with $n_a >0$. For roots $a \in  \Delta \backslash Z$  we have $a(z_n) \rightarrow \infty$. Hence the second term $(\lambda_0( \Delta) - \lambda)(z_n - u_n)$ goes to $\infty$, and for $n$ big enough we have  $(\lambda_0( \Delta') - \lambda)(z_n - u_n) > 0$. 

Therefore we have $z_n - u_n \in F'$ for $n$ big enough, which implies that almost all $z_n$ are contained in $U + F'$. Thus almost all $x_n$ are contained in $\Gamma_{F'}(U)$.

Hence we find that $x_n$ converges to $x$, which completes the proof of compactness.\hfill$\Box$

The apartment $A$ is equipped with an action
\[\nu: N \longrightarrow \mbox{\rm Aff}(A),\]
 of $N$, where $\mbox{\rm Aff}(A)$ denotes the affine-linear automorphisms of $A$. 
Let $v_K$ be the valuation on $K$, normalized
so that a prime element has valuation one. 
An element $z \in Z$ acts by translation with the element 
$\nu(z)$ in $A$ satisfying $\chi(\nu(z)) = -v_K(\chi(z))$ for all $\chi \in X^\ast(\bz)$. 
Since $X^\ast(\bz) \subset X^\ast(\bs)$ is of finite index, this condition defines indeed a point in $A$.
For $n \in N$ we denote by $\ov{n}$ the induced element in the Weyl group $W = N/Z$. 
The linear part of every affine automorphism $\nu(n)$ is given by the natural action of $\ov{n} \in W$ on $A$.

Now we want to extend $\nu$ to $\ov{A}$.
First of all, note that every $w \in W$ acts on the set of faces $\Sigma$,
since $w$ maps $F = \fyd$ to $w(F) = F^{w(\Delta)}_{w(Y)}$.
Hence for every $n \in N$ with linear part $w = \ov{n} \in W$ 
the affine automorphism $\nu(n)$ of $A$ induces an affine map 
\[A_F = A / \langle F \rangle \longrightarrow A_{w(F)} = A / \langle w(F)\rangle.\]
If we put all these affine maps on boundary components together,
we get a map $\nu(n): \ov{A} \rightarrow  \ov{A}$.
This defines an action of $N$ on $\ov{A}$. 

\begin{lem}
The action $N \times \ov{A} \rightarrow \ov{A}$ defined by $\nu$ 
is continuous.
\end{lem}
{\bf Proof: } For every $n \in N$ we have $\nu(n)(\Gamma_F(U)) = \Gamma_{\ov{n}(F)}(\nu(n)(U))$, hence $\nu(n)$ is a continuous
map on $A$. By construction, the open subgroup $\{z \in Z: \chi(z) \in \mathfrak{O}_K^\times \mbox{ for all } \chi \in X^\ast(\bz)\}$
of $N$ acts identically on $A$, and hence on $\ov{A}$. Therefore the action of $N$ on $\ov{A}$ is continuous.\hfill$\Box$

\section{The stabilizer groups $P_x$}
Let us first recall the construction of the Bruhat-Tits building $X = X(\bg)$. For every root $a \in \Phi$ denote by $\bu_a$ 
the root group (see \cite{bo}, 21.9, where the notation is $\bu_{(a)}$), and write $U_a$ for its group of $K$-rational points.
Bruhat and Tits define an exhaustive and separating filtration on $U_a$ by  non-trivial open compact subgroups $U_{a,r}$ for $r \in \R$, see \cite{brti},  (6.1.3) and \cite{brti2}, (5.2.2). This filtration decreases with increasing $r$,
i.e. $U_{a,r} \supset U_{a,s}$ if $r \leq s$, and jumps occur only in a discrete subset of $\R$.
We put $U_{a, \infty} = \{1\}$ and $U_{a, -\infty} = U_a$.

For every $x \in A$ and every root $a \in \Phired$ define $U_{a,x} = U_{a, -a(x)}$,
and let $U_x$ be the subgroup of $G$ which is generated by all $U_{a,x}$ for $a \in \Phired$.
Besides, put $N_x = \{ n \in N: \nu(n) x = x\}$ and $P_x = N_x U_x$. Then the building is defined as follows: 
\[ X = X(\bg) = G \times A /\sim,\]
where the equivalence relation $\sim$ is defined by
\begin{eqnarray*}
(g,x) \sim (h,y) & \mbox{if and only if  there exists an element } n \in N \\
~ & \mbox{such that } \nu(n)x= y \mbox{ and } g\inv h n \in P_x.
\end{eqnarray*}
There is a natural action of $G$ on $X(\bg)$ via left multiplication on the first factor. 
Besides, $X(\bg)$ carries a  topology induced by the topology on $G \times A$. Note
that in this product the first factor is of $p$-adic nature, whereas the second one is a real vector space.
The apartment $A$ can be embedded in $X(\bg)$ by mapping $a \in A$ to the class of $(1,a)$.
The $G$-action on $X$ continues the $N$-action on $A$, so that we simply write $nx$ instead of $\nu(n) x$ if $x \in A$ and $n \in N$.  

Recall that we fix a geometrically irreducible faithful $K$-rational linear representation $\rho$ of $\bg$. As shown in the previous section,
$\rho$ defines a compactification $\ov{A}$ of the apartment $A$.  For every subset $X$ of $A$ we denote by $\ov{X}$ the closure of $X$ in $\ov{A}$.

\begin{defi} Let $\Omega$ be a subset of $\ov{A}$. For every root $a \in \Phired$ we put
\begin{eqnarray*}
 f_\Omega(a) & =  \inf\{t: \Omega \subseteq\overline{\{ z \in A: a(z) \geq -t\}} \} \\
~ & = -\sup \{ t: \Omega \subseteq \ov{\{ z \in A: a(z) \geq t \} } \}.
\end{eqnarray*}
\end{defi}

Here $\inf \emptyset = \sup \R = \infty$ and $\inf \R = \sup \emptyset = - \infty$. If $\Omega = \{x\}$, we write $f_x(a) = f_{\{x\}}(a)$. Note that 
\[ f_x(a) = -a(x) \quad \mbox{ for all }x \in A,\]
and that 
\[f_{\om_1}(a) \leq f_{\om_2}(a), \quad \mbox{if }\om_1 \subseteq \om_2.\]
For all $a \in \Phired$ with $2a \in \Phi$ we put $f_\Omega(2a) = 2 f_\Omega(a)$ with obvious rules for $\pm \infty$. Then
$f_\Omega$ is a function $f_\Omega: \Phi \rightarrow \R \cup \{\infty, - \infty\}$. 

If  $F = \fyd$ is a face, we identify as in the last
section $\langle Y \rangle$ with the dual space of $A_F$, hence we can apply elements
of $\langle Y \rangle$ to points in $A_F$. It is easy to check that  for every $x \in A_F$ and
every $a \in \langle Y \rangle$ we have
$f_x(a) = -a(x)$. 

\begin{defi}
For every $\Omega \subset \ov{A}$ and $a \in \Phi$ put $U_{a, \Omega} = U_{a, f_\Omega(a)}$. Let $U_\Omega$ denote the subgroup of $G$ which is generated by all groups $U_{a, \Omega}$ for $a \in \Phired$.
\end{defi}

Note that for $\Omega = \{x\}$ with $x \in A$ we get the groups $U_{a,x}$ and $U_x$ from above. 
We will now investigate the extreme cases $U_{a, \Omega} = 1$ and $U_{a, \Omega} = U_a$.

\begin{prop} Let $\Omega$ be a non-empty subset of $\ov{A}$ and $a \in \Phired$.

i) $U_{a, \Omega} = U_a$ if and only if  for all faces $F = \fyd$ such that $\Omega \cap A_F \neq \emptyset$ the following condition holds: $a \notin \langle Y \rangle$ and for all $K$-weights $\lambda$ with $[\lambda_0( \Delta) - \lambda] \subset Y$ and all positive integers $l$ the element $\lambda + l a$ is not a $K$-weight.

ii) $U_{a, \Omega} = \{1\}$ if and only if there exists a face $F = \fyd$ with $\Omega \cap A_F \neq \emptyset$
such that one of the following conditions hold: Either $a \in \langle Y \rangle$ and $a(\Omega \cap A_F)$ is unbounded from below or 
$a \notin \langle Y \rangle$ and there exists a $K$-weight $\lambda$ with $[\lambda_0( \Delta) - \lambda] \subset Y$ and a positive integer $l$ such that $\lambda + la$ is a $K$-weight.
\end{prop}

{\bf Proof: } 
i) By definition $U_{a, \Omega} = U_a$, if and only if $f_{\Omega}(a) = - \infty$, which means
that $\Omega \subset \ov{\{ a \geq t\} }$ for all $t \in \R$. 
Assume that $f_\Omega(a) = - \infty$, and let $x \in \Omega \cap A_F$ for the face $F = \fyd$.
Since $f_x(a) = - \infty$, there exists a sequence of points $x_n \in A$ converging towards $x$ such that $a(x_n) \rightarrow \infty$. On 
the other hand, 
the sequence $b(x_n)$ is bounded for all $b \in \langle Y \rangle$.
This implies $a \notin \langle Y \rangle$. 

Suppose $\lambda$ is a $K$-weight such that $[\lambda_0( \Delta) - \lambda] \subset Y$. Let $l \geq 1$. Then $(\lambda_0(\Delta) - \lambda -la)(x_n) \rightarrow -\infty$.
On the other hand, for every $K$-weight $\mu$ with $[\lambda_0(\Delta) - \mu] \not\subset Y$ we have $(\lambda_0(\Delta)-\mu)(x_n)\rightarrow \infty$. Hence $\lambda + la$ is not a $K$-weight. 

Let us now assume that $x$ is a point  in $A_F$ for $F = \fyd$ such that $a \notin \langle Y \rangle$ and such that for all $K$-weights $\lambda$ with $[\lambda_0 (\Delta) - \lambda] \subset Y$ and all $l \geq 1$  the element $\lambda + la$ is not a $K$-weight. We have to show $f_x(a) = - \infty$, i.e. $x \in \ov{\{ a \geq t\}}$ for all $t \in \R$. We denote by $\Phi^+$ respectively $\Phi^-$ the positive respectively negative roots with respect to $\Delta$. 

If $a \in \Phi^+$, let $x_n$ be a sequence of points in $A$ such that 
\begin{eqnarray*}
b(x_n) = b(x),  & \mbox{ if } & b \in Y \mbox{ and}\\
b(x_n) \rightarrow \infty, & \mbox{ if }& b \in \Delta \backslash Y
\end{eqnarray*}
Then $x_n$ converges to $x$ in $\ov{A}$, and $a(x_n) \rightarrow \infty$, which implies $f_x(a) = - \infty$. 

If $a \in \Phi^-$, we have $-a = \sum_{b \in  \Delta} n_b b$ with $n_b \geq 0$. Put $[-a] = \{ b \in  \Delta: n_b > 0\}$. 

{\bf Claim: }Let $a \in \Phi^-$ such that  $Y \cup [-a]$ is admissible. Then there exists a $K$-weight $\lambda$ with $[\hiwei - \lambda] \subset Y$ and some positive integer $l$ such that $\lambda + l a$ is a $K$-weight. 

Let us prove this claim. 
The weights of $\rho_\kbar$ are stable under $\lkbar W$. By \cite{boti}, 5.5, for every $w \in W$ there is an element $w' \in \lkbar W $ such that  $w'$ restricted to $\bs_\kbar$ is equal to $w$. Hence the $K$-weights are stable under the action of $W$. If $w = r_a \in W$ is the reflection corresponding to the root $a$ and $\lambda$ is a $K$-weight, we have $r_a (\lambda) = \lambda -2 (\lambda, a) / (a,a) a$, where $(\, , \,)$ denotes a $W$-invariant scalar product. Hence, if $\lambda = \lambda_0(\Delta)$ and $b \in \Delta$, then $(\lambda_0(\Delta), b) \geq 0$. 

Assume that $Y \cup [-a]$ is admissible, i.e. that $Y \cup [-a] \cup \{\lambda_0( \Delta)\}$ is connected. Since $a \in \Phi^-$, we have $(\lambda_0(\Delta),a) \leq 0$.
If $(\lambda_0(\Delta),a) <0$, then $r_a (\lambda_0(\Delta))$ is a $K$-weight of the form $\lambda_0(\Delta) + la$ for some $l \geq 1$. 
If $(\lambda_0(\Delta),a) =0$,  we find $a_1, \ldots, a_r \in Y \backslash [-a]$ with $\{a_1, \ldots, a_r\}$ admissible and  $(a_r, b) \neq 0$ for some $b \in [-a]$.  
Hence there exists a $K$-weight $\lambda = \hiwei - \sum_{i = 1}^r  m_i a_i$ with $m_i >0$ by \cite{boti}, 12.16. 
Since all $a_i$ and $b$ belong to $ \Delta$, we have $(a_i, b) \leq 0$  by \cite{bou}, VI, 1.3. In particular, $(a_r, b) < 0$. Hence
\[(\lambda, a) = (\hiwei- \sum_i m_i a_i, -\sum_{c \in [-a]} n_c c)
= -\sum_c n_c (\hiwei, c) + \sum_i \sum_c m_i n_c (a_i,c).\]
Since the first sum is $\leq 0$ and all $(a_i,c)$ are $\leq 0$ with $(a_r,b) <0$, we find that
$(\lambda, a) < 0$. Hence $r_a(\lambda )$ is a $K$-weight of the form $\lambda +l a $ for some $l \geq1$, which proves the claim.

Now we finish the proof of i) by showing that $f_x(a) = -\infty$ for $a \in \Phi^-$. By the claim, 
our assumption  implies that $Y \cup [-a]$ is not admissible. Note that $Y \cup \{\hiwei\}$ is connected by definition, and $[-a]$ is connected by \cite{bou}, VI, 1.6. Therefore $Y$ and $[-a]$ are disjoint.
Let  $\lambda$ be a $K$-weight such that $[\hiwei - \lambda] \not\subset Y$. Then 
$[\hiwei - \lambda]$ is admissible, i.e. $[\hiwei - \lambda] \cup \{\hiwei\}$ is connected.
Hence $[\hiwei - \lambda]$ cannot be contained in $Y \cup [-a]$. 
We write
\[ \hiwei - \lambda = \sum_{b \in Y} m_b b + \sum_{b \in [-a] } m_b b + \sum_{b \notin (Y \cup [-a])} m_b b\]
with $m_b \geq 0$. The
last sum is not zero, i.e. for at least one $b \notin (Y \cup [-a])$ we have $m_b >0$. 
Since there are only finitely many $K$-weights, there is a natural number $N$ satisfying
\[N \sum_{b \notin (Y \cup [-a])} m_b \geq \sum_{b \in [-a] } m_b +1\]
for all $\lambda$ such that $[\hiwei - \lambda] \not\subset Y$. 

Now let $x_n$ be a sequence of points in $A$ satisfying
\begin{eqnarray*}
b(x_n) = b(x) & \mbox{ for } & b \in Y\\
b(x_n) = Nn & \mbox{ for } & b \in  \Delta \backslash (Y \cup [-a])\\
b(x_n) = -n & \mbox{ for  } &  b \in [-a] .
\end{eqnarray*}
For all $K$-weights $\lambda$ with $[\hiwei - \lambda] \not\subset Y$ we have 
\[(\hiwei- \lambda)(x_n) = \sum_{b \in Y} m_b b(x) + n \left(- \sum_{b \in [-a]} m_b  + N \sum_{b \notin (Y \cup [-a])} m_b\right),\]
hence $(\hiwei - \lambda) (x_n) \rightarrow \infty$. Therefore $x_n$ converges towards $x$ in $\ov{A}$.
Since $a(x_n) \rightarrow \infty$, we conclude that $f_x(a) = - \infty$.

ii) By definition, $U_{a, \Omega} = 1$ is equivalent to $f_\Omega(a) = \infty$, which means
that $\Omega$ is not contained in any set of the form $\ov{\{a \geq t\}}$. 

Assume that $f_\Omega(a) = \infty$. Then there is some face $F = \fyd$ with $\Omega \cap A_F \neq \emptyset$ and $f_{\Omega \cap A_F} = \infty$.  If $a$ is contained in $\langle Y \rangle$, we have
$f_x(a) = -a(x)$ for all $x \in \Omega \cap A_F$, hence $a(\Omega \cap A_F)$ must be unbounded from below.  If, on the other hand, $a$ is not contained in $\langle Y \rangle$, then there exists a $K$-weight $\lambda$ as in the statement, since otherwise we could use i) to deduce $f_{\Omega \cap A_F} = - \infty$.

Let us prove the other direction. If $a \in \langle Y \rangle$ and $a(\Omega \cap A_F)$ is unbounded from below, then
obviously $f_{\Omega \cap A_F}(a) = \infty$, hence $f_\Omega(a) = \infty$. If $a \notin \langle Y \rangle$, assume that there exists a $K$-weight $\lambda$ and some $l\geq 1$ with $[\hiwei - \lambda] \subset Y$ such that $\lambda +l a$ is a $K$-weight. Let $x_n \in A$ be any sequence converging to some point $x \in \Omega \cap A_F$. Then
$(\hiwei - \lambda)(x_n)$ is bounded, whereas $(\hiwei - \lambda -la)(x_n) \rightarrow \infty$. 
This implies $a(x_n) \rightarrow -\infty$. Hence $x$ is not contained in any set of the form $\ov{\{a(z) \geq t\}}$, which means $f_x(a) = \infty$ and also $f_\Omega(a) = \infty$. \hfill$\Box$

Recall that we denote the quotient map $N \rightarrow W$ by $n \mapsto \ov{n}$. For $n \in N$  and the root group $U_a$ we have $n U_a n\inv = U_{\ov{n} (a)}$.

\begin{lem}
Let $\Omega$ be a non-empty subset of the compactified apartment $\ov{A}$. For every $n\in N$ and $a \in \Phi$  we have $n U_{a,\Omega} n\inv = U_{\ov{n}(a), n(\Omega)}$. Hence 
$n U_\Omega n\inv = U_{n (\Omega)}$.
\end{lem}

{\bf Proof: } 
By \cite{brti}, 6.2.10, the claim holds for subsets $\Omega$ of $A$. 
Since $n$ acts on $A$ by affine automorphisms,  for all $z \in A$ we have $n(z) = n(0) + \ov{n}(z)$. Hence 
$a(n(z)) = a(n(0)) + a(\ov{n}(z))$.
Therefore
\[n\{ z\in A: a (z) \geq s\} = \{z \in A: \ov{n}(a)(z) \geq \ov{n}(a)(n(0)) +s\}.\] 
Since
$n$ acts as a homeomorphism on $A$, we also have
\[ n\ov{\{ z \in A:  a(z) \geq s\}} = \ov{\{z \in A: \ov{n}(a)(z) \geq \ov{n}(a)(n(0)) +s\}}.\]
Hence $f_\Omega(a) = \pm \infty$ if and only if $f_{n(\Omega)} ( \ov{n}(a) ) = \pm \infty$. Therefore we can assume $f_\Omega (a) \in \R$. Then there exists a face $F= \fyd$ such that $f_{\Omega \cap A_F}( a) = f_{\Omega}(a)$. By Proposition 3.3, we deduce $a \in \langle Y \rangle$. Put $\Omega' = r_F\inv(\Omega \cap A_F) \subset A$. Then $f_{\Omega'}(a) = f_{\Omega \cap A_F}(a) = f_\Omega(a)$ and $f_{n(\Omega')}(\ov{n}(a))= f_{n(\Omega)\cap A_{\ov{n} F}}(\ov{n}(a)) = f_{n(\Omega)}(\ov{n}(a))$. Hence our claim follows from the corresponding result for subsets of $A$.\hfill$\Box$

\begin{lem}
Let $x_k$ be a sequence of points in $\ov{A}$, which converges to $x \in \ov{A}$. Let $u_k \in U_{a,x_k}$ be a sequence of elements, converging to some $u$ in $U_a$. Then $u$ is contained in  $U_{a,x}$.
\end{lem}
{\bf Proof: }If $f_x(a) = - \infty$, there is nothing to prove, since  $U_{a,x} = U_a$.
 
Assume that $f_x(a) = \infty$. Then any set $\{\ov{z \in A: a(z) \geq s}\}$ 
 contains only finitely many elements $x_k$. Hence the sequence $f_{x_k}(a)$ goes to $\infty$, so that $u \in \cap_r U_{a,r} = \{1\}$. 

If $f_x(a) \in \R$, we have $x \in A_F$ for some face $F = \fyd$ with $a \in \langle Y \rangle$ by Proposition 3.3. Hence $f_x(a) = -a(x)$. Then almost all $x_k$ are contained in boundary components $A_{F'}$ for faces $F' = F_{Y'}^{\Delta'}$ such that $\langle Y \rangle \subset \langle Y' \rangle$. Hence $f_{x_k}(a) = -a(x_k)$, which implies $f_{x_k}(a) \rightarrow f_x(a)$. Assume that $U_{a,r}$ is the last filtration step containing $u$, i.e. $u \in U_{a,r}$, but $u \notin U_{a, r+ \varepsilon}$ for all $\varepsilon > 0$. Then for all $\varepsilon > 0$ almost all $u_k$ are not contained in the closed subgroup $U_{a, r + \varepsilon}$. Hence for almost all $k$ we have $r + \varepsilon > f_{x_k}(a)$, which implies $r + \varepsilon \geq f_{x}(a)$. Therefore  $r \geq f_{x}(a)$, so that $u \in U_{a,r} \subset U_{a,x}$.\hfill$\Box$

For the next result we choose a Weyl chamber of the root system $\Phi$ and denote the corresponding
positive respectively negative roots by $\Phi^+$ respectively $\Phi^-$. Besides, we put $\Phiredp = \Phired \cap \Phi^+$ and $\Phiredm = \Phired \cap \Phi^-$. 
Denote by $\bu^+$ respectively $\bu^-$ the corresponding subgroups of $\bu$, see \cite{bo}, 21.9 and
by $U^+$ and $U^-$ their $K$-rational points. Then $U^\pm$ is directly spanned by the $U_a$ for $a \in \Phiredpm$ in any order.

For every subset $\om$ of $\ov{A}$ we put 
\[U_\om^+ = U_\om \cap U^+ \mbox{ and } U_\om^- = U_\om \cap U^-.\]

\begin{thm}
i) The product map induces a bijection \[\prod_{a \in \Phiredpm} U_{a, \Omega} \longrightarrow U_{\Omega}^\pm,\] where the product on the left hand side may be taken in arbitrary order.

ii) $U_a \cap U_\om = U_{a, \om}$ for all $a \in \Phired$.

iii) $U_\om = U_\om^- U_\om^+ (N \cap U_\om)$.
\end{thm}
{\bf Proof: }  We adapt the arguments in the proof of \cite{brti}, 6.4.9, to our situation.
Note that for $a \in \Phi$ with $2a \in \Phi$ we have $U_{2a,\om} \subset U_{a,\om}$ by \cite{brti}, 6.2.1.

Let $a$ and $b$ be non-proportional  roots in $\Phired$.
We first claim that the commutator $(U_{a,\om}, U_{b,\om})$ is contained in the group generated by all
$U_{ ma+nb, \om}$ for $m, n \geq 1$ such that $ma + nb \in \Phi$.

If $f_{\om}(a) = \infty$ or $f_\om(b) = \infty$, the commutator $(U_{a,\om}, U_{b,\om})$ is trivial,
so that  this claim holds. 
Hence we can assume that for all faces $F$ with $\om \cap A_F \neq \emptyset$ we have  $f_{\om \cap A_F}(a) \in \R \cup \{ - \infty\}$ and  $f_{\om \cap A_F}(b) \in \R \cup\{ -\infty\}$.

First we treat the case that for all faces $F$ with $\om \cap A_F \neq \emptyset$ we have  $f_{\om \cap A_F}(a) = - \infty$ or  $f_{\om \cap A_F}(b) = -\infty$. Let $ma+nb$ with $m,n \geq 1$ be a root. We want to show that in this case $f_{\om \cap A_F}(ma+nb) = - \infty$.

First we treat the case  $f_{\om \cap A_F }(a) = -\infty$ and $f_{\om \cap A_F} (b) = -\infty$. 
By Proposition 3.3, we find $a \notin \langle Y \rangle$ and $b \notin \langle Y \rangle$. If $a$ and $b$ are both positive or both negative with respect to the ordering defined by $\Delta$, i.e. $a, b \in \Phi^+(\Delta)$ or $a,b \in \Phi^-(\Delta)$, then obviously $ma+nb \notin \langle Y \rangle$. If $ a \in \Phi^+(\Delta)$ and $b \in \Phi^-(\Delta)$, we also claim that $ma+nb \notin \langle Y \rangle$. Let us again denote for every positive root $c = \sum_{d \in \Delta} n_d d$ its support by $[c] = \{d \in \Delta: n_d > 0\}$. By \cite{bou}, VI, 6.1, the set $[c]$ is connected. Now assume that $ma+nb \in \langle Y \rangle$ and $ma+nb \in \Phi^+(\Delta)$. Then  $[ma+nb] \subset Y$. Either $[-b] \cap [ma+nb] \neq \emptyset$ or $[-b] \cup [ma+nb] = [a]$. Since $[a]$ and $[-b]$ are connected and $[ma+nb]$ is contained in the admissible set $Y$,  we deduce in both cases that $Y \cup [-b]$ is admissible. This contradicts  the claim in Proposition 3.3, which says that the condition $f_{\om \cap A_F}(b) = -\infty$ implies that $Y \cup [-b]$ is not admissible.  The other cases can be treated similarly. 
Hence we find $ma+ nb \notin \langle Y \rangle$. 

Therefore by Proposition 3.3 we have $f_{\om \cap A_F}(ma+nb) = \infty$ or $-\infty$.  If $ma+nb$ lies in $\Phi^+(\Delta)$, we also have $f_{\om \cap A_F}(ma+nb) = - \infty$. Let us   assume that $ma+nb \in \Phi^-(\Delta)$. If both $a$ and $b$ lie in $\Phi^-(\Delta)$, we have
$[-ma-nb] = [-a] \cup [-b]$. Since neither $Y \cup [-a]$ nor $Y \cup [-b]$ is admissible, the set $Y \cup [-ma-nb]$ is also not admissible. By Proposition 3.3 we find $f_{\om \cap A_F}(ma+nb) = - \infty$. 
If $a$ lies in $\Phi^-(\Delta)$ and $b$ lies in $\Phi^+(\Delta)$ (and similarly, if it is the other way round),
we find $[-ma-nb] \subset [-a]$. Since $Y \cup [-a]$ is not admissible and $[-a]$ is connected, the subset $[-ma-nb]$ is not admissible. Hence again $f_{\om \cap A_F}(ma+nb) = -\infty$. 

Now we consider a face $F = \fyd$ with   $\om \cap A_F \neq \emptyset$ such that $f_{\om \cap A_F}(a) = -\infty$ and $f_{\om \cap A_F}(b) \in \R$. By  Proposition 3.3, $a \notin \langle Y \rangle$ and $b \in \langle Y \rangle$, so that $ma+nb \notin \langle Y \rangle$.
 Again, if $ma+nb$  lies in $\Phi^+(\Delta)$, then Proposition 3.3 implies $f_{\om \cap A_F}(ma+nb) = - \infty$.  If both $a$ and $b$ lie in $\Phi^-(\Delta)$, we have
$[-ma-nb] = [-a] \cup [-b]$. Since   $Y \cup [-ma-nb]= Y \cup [-a]$ is not admissible, we use Proposition 3.3 to deduce $f_{\om \cap A_F}(ma+nb) = - \infty$. 
If $a$ lies in $\Phi^-(\Delta)$ and $b$ lies in $\Phi^+(\Delta)$, again $Y \cup [-ma-nb] = Y \cup [-a]$
is not admissible, so that $f_{\om \cap A_F}(ma+nb) = -\infty$. The case $a \in \Phi^+(\Delta) $ and $b \in \Phi^-(\Delta)$ occurs only for $ma+nb \in \Phi^+(\Delta)$.  
Therefore we have seen that $f_{\om \cap A_F}(a) = - \infty$ and $f_{\om \cap A_F}(b) \in \R$ imply
$f_{\om \cap A_F}(ma+nb) = -\infty$. The same argument with reversed roles shows that $f_{\om \cap A_F}(a) \in \R$ and $f_{\om \cap A_F}(b) = -\infty$ implies
$f_{\om \cap A_F}(ma+nb) = -\infty$.  

Hence we deduce that $f_\om(ma+nb) = - \infty$, if for all faces $F$ with $\om \cap A_F \neq \emptyset$ at least one of the numbers $f_{\om \cap A_F}(a)$, $f_{\om \cap A_F}(b)$ is $-\infty$. In this case 
$U_{ma+nb, \om} = U_{ma+nb}$, and our claim follows from the analogous statement for the groups $U_a$, see  \cite{bo}, 14.5 and 21.9.

It remains to treat the case that  for some $F = \fyd$ with $\om \cap A_F \neq \emptyset$, we have $f_{\om \cap A_F}(a) \in \R $ and $f_{\om \cap A_F}(b) \in \R$.  By Proposition 3.3, we have $a, b \in \langle Y \rangle$, and also $ma+ nb \in \langle Y \rangle$, whenever $ma+nb$ with $m,n  \geq 1$ is a root. Let  $r,s$ be real numbers with $r > f_{\om \cap A_F}(a)$ and $s > f_{\om \cap A_F}(b)$. Then $\om \cap A_F$ is contained in 
\[A_F \cap \ov{\{z \in A: a(z) \geq -r\}} = \{y \in A_F: a(y) \geq -r\}\]
and in 
\[A_F \cap \ov{\{z \in A: b(z) \geq -s\}} = \{y \in A_F: b(y) \geq -s\}.\]
Hence $\om \cap A_F \subset \ov{\{z \in A: ma+nb(z) \geq -mr-ns\}}$.
Therefore $f_{\om \cap A_F}(ma+nb) 
\leq m f_{\om \cap A_F}(a) + n f_{\om \cap A_F}(b)$.  
Since $f_\om(a)$ is the maximum of all $ f_{\om \cap A_F}(a)$ for $\om \cap A_F \neq \emptyset$, we find 
$f_\om(ma+nb) \leq  m f_\om (a)+ n f_\om (b)$.  Then the claim  follows from  \cite{brti}, 6.2.1, (V3).

Let $H^\pm$ be the subgroup of $U_\om$ generated by all $U_{a, \om}$ for $a \in \Phiredpm$.
By \cite{brti}, 6.1.6  the product map induces a bijection 
\[\prod_{a \in \Phiredpm} U_{a, \Omega} \rightarrow H^\pm. \]
Let $U_\om^{(a)}$ be the subgroup of $G$ generated by $U_{a, \om}$ and by $U_{-a, \om}$, and put
$N_\om^{(a)} = N \cap U_\om^{(a)}$. Let $Y$ be the subgroup of $U_\om$ generated by all $N_\om^{(a)}$ for $ a \in \Phired$. 

Now we show that the product $H^- H^+ Y$ is independent of the choice of the Weyl chamber $D$ corresponding to the chosen order on $\Phi$. If suffices to show that the product remains unchanged when $D$ is replaced by $r_a(D)$, where $r_a$ is the reflection associated to the root $a$ in the basis of $\Phi$ given by $D$. Denote by $H'_a$ (respectively $H'_{-a}$) the product of all $U_{b,\om}$ for $b \neq a$ in $\Phiredp$ (respectively $b \neq -a$ in $\Phiredm$).
By \cite{brti}, 6.1.6,  $H'_a$ and $H'_{-a}$ are groups. The commutator relation proved at the 
beginning implies that $H'_a$ and $H'_{-a}$ are normalized by  both $U_{-a, \om}$ and $U_{a, \om}$.

Now we show that $U_{ -a, \om} U_{a, \om} \subset U_{a, \om} U_{-a, \om} N_\om^{(a)}$.
If $f_\om(a) = \infty$ or $f_\om(-a) = \infty$, this is clear. Hence we can assume that neither $f_\om(a)$ nor $f_\om(-a)$ is $\infty$. 
Then $f_\om(a) \in \R$ if and only if $f_\om(-a) \in \R$, and in this case $f_\om(a) + f_\om(-a) \geq 0$. Hence the function $f_\om : \{a,-a\} \rightarrow \R$
can be continued to a concave function $\Phi \rightarrow \R$. By \cite{brti}, 6.4.7, this implies that
$U_\om^{(a)} = U_{a, \om} U_{-a, \om} N_\om^{(a)}$. Hence the desired inclusion holds.
It remains  to consider the case $f_\om(a) = f_\om(-a) =  - \infty$, i.e. $U_{a,\om} = U_a$ and $U_{-a, \om} = U_{-a}$. For all $u \in U_{-a} \backslash \{1\}$ there exists an element $m(u) \in N$ with $u \in U_a m(u) U_a$ by  \cite{brti}, 6.1.2 and 6.1.3.
Then $m(u) \in N_\om^{(a)}$, and there are $v_1, v_2 \in U_{a}$ with $u =  v_1 m(u) v_2$. Hence for all $v \in U_{a}$
we have 
\[u v = v_1 m(u) v_2 v = v_1 \, (m(u) \, v_2 v \, m(u)\inv ) \, m(u).\]
By \cite{brti}, 6.1.2, $m(u) \, v_2 v \, m(u)\inv$ lies in $U_{-a}$, which implies $uv \in U_a U_{-a} N_\om^{(a)}$. Hence also in this case $U_{-a, \om} U_{a, \om} \subset U_{a, \om} U_{-a, \om} N_\om^{(a)}$.

Note that this inclusion implies $U_{-a, \om} U_{a, \om} Y = U_{a, \om} U_{-a, \om} Y$.

Therefore we have 
\begin{eqnarray*}
H^- H^+ Y & = & H'_{-a} U_{-a, \om} U_{a, \om} H'_{a} Y = H'_{-a} H'_a U_{-a, \om} U_{a, \om} Y \\
 & = & H'_{-a}H'_a U_{a, \om} U_{-a, \om} Y = H'_{-a} U_{a, \om} H'_a U_{-a, \om} Y, 
\end{eqnarray*}
which shows that $H^- H^+ Y$ remains unchanged when we change the order on $\Phi$ from $D$ to $r_a(D)$. 

It follows that $H^- H^+ Y$ is stable under multiplication from the left with all $U_{a, \om}$ for $a \in \Phired$. 
Therefore $U_{\om} = H^- H^+ Y$.

Since $U^+ \cap U^- = \{1\}$ and $N \cap U^- U^+ = \{1\}$  by \cite{boti}, 5.15, we find $U_\om^- = U_\om \cap U^- = 
 H^- H^+ Y \cap U^- = H^-$. Similarly, $U_\om^+ = H^+$, which shows i). For all $a \in \Phiredpm$ we have
$U_a \cap U_\om = U_a \cap H^\pm= U_{a, \om}$, which proves ii).
The only point that  remains to be checked  for iii) is $N \cap U_\om = Y$. Obviously, $Y \subset N \cap U_\om$. 
Since $N \cap U^- U^+ = \{1\}$, we get $N \cap U_\om = N \cap  H^- H^+ Y  \subset Y$, which finishes the proof.\hfill$\Box$

Now we denote for any non-empty $\Omega \subset \ov{A}$  by $N_\Omega = \{ n \in N: nx=x \mbox{ for all } x \in \Omega\}$ the stabilizer of $\Omega$ in $N$ with respect to the action $N \times \ov{A} \rightarrow A$ defined in the previous section. We write $N_x=N_{\{x\}}$ for all $x \in \ov{A}$. 
\begin{cor}
We have $N \cap U_\Omega \subset N_\Omega$. 
\end{cor}
{\bf Proof: }
Since for all $x \in \Omega$ we have $U_\Omega \subset U_x$, it suffices to show the claim for one point sets $\Omega = \{x\}$. 
We have seen in the proof of the theorem that $N \cap U_x$ is the subgroup of $U_x$ generated by all $N \cap U_x^{(a)}$, where $U_x^{(a)}$ is the subgroup generated by $U_{a,x}$ and $U_{-a,x}$.
Hence we have to show that for all $a \in \Phired$ the group $N \cap U_x^{(a)}$ stabilizes $x$.
Let $A_F$ for $F = \fyd$ denote the stratum of $\ov{A}$ containing $x$. 

If $f_x(a) = \infty$ or $f_x(-a) = \infty$, then $U_{a,x}$ or $U_{-a,x}$ is trivial. Hence $N \cap U_x^{(a)}$ is trivial, and our claim is true.
We can therefore assume that $f_x(a)$ and $f_x(-a)$ are contained in $\R \cup \{-\infty\}$. 

If $f_x(a)$ is contained in $\R$, then by Proposition 3.3 also $f_x(-a)$ is contained in $\R$ and $a \in \langle Y \rangle$. For every $z \in A$ with $a(z) = a(x)$ we find $U_{a,x} = U_{a,z}$ and $U_{-a,x} = U_{-a,z}$. Hence every $n \in N \cap U_x^{(a)}$ is contained in $N \cap U_z$. Therefore $n$ stabilizes the point $z$, which is contained in the Bruhat-Tits building. Since we can find a sequence $z_n$ in $A$ converging to $x$ with  $a(z_n) = a(x)$ and since $N$ acts continuously on $\ov{A}$, it follows that $n$ stabilizes $x$. 

It remains to consider the case $f_x(a) = -\infty$ and $f_x(-a) = - \infty$. We can exchange $a$ and $-a$ and may therefore assume that $a \in \Phiredp$. Then the claim in Proposition 3.3 implies that $[a] \cup Y$ is not admissible. By \cite{brti}, 6.1.2 and 6.1.3 the group generated by $U_a$ and $U_{-a}$ (and hence also $U_x^{(a)}$) is contained in the union $Z U_a \cup U_a m Z U_a$, where  $m Z$ is the coset in $N$ given by some element $m = m(u)$ for $u \in U_{-a} \backslash \{1\}$, i.e. $m(u)$ is the unique element in in $U_{a} u U_{a} \cap N$. Then the image of $m(u)$ in the Weyl group is  the reflection corresponding to $a$. 
Fix some $n \in N \cap U_x^{(a)}$.  

Let us first consider the case that ${n}$ is contained in $Z U_a$, hence in $Z$. Let $b$ be a root in $Y$. Then $b$ is perpendicular to $[a]$, hence $U_a$ and $U_{-a}$ are contained in the  centralizer of $U_b$ by \cite{bo}, 14.5 and 21.9. Therefore $n$ centralizes $U_b$. An analogous argument shows that $n$ centralizes $U_{-b}$. Now let  $u$ be an element of $U_{b} \backslash \{1\}$. By \cite{ti}, 1.4, $m(u)$  acts on $A$ as an affine reflection along some affine hyperplane $\{ z \in A: b(z) + r = 0\}$. Now $m(n\inv u n) = m(u)$, and on the other hand $m(n\inv u n) = n\inv m(u) n$ acts on $A$ as the reflection corresponding to the affine hyperplane $\{z \in A: b(z) + r -b(\nu(n)) = 0\}$. Hence $b(\nu(n)) = 0$. We consider a sequence of points $z_k \in A$ satisfying
$b(z_k) = b(x)$ for all $b \in Y$ and $b(z_k) \rightarrow \infty$ for all $b \in \Delta \backslash Y$. Then $z_k$ converges towards $x$ and $n z_k$ converges towards  $x$. As $Z$ acts continuously on $\ov{A}$ by Lemma 2.5, we deduce $nx = x$. 

Now we consider the case $n \in U_a m(u) h U_a$ for some $h \in Z$. Then $n (m(u) h)\inv \in U_a m(u) h U_a (m(u) h)\inv = U_a U_{-a}$ by \cite{brti}, 6.1.2.
Since $N \cap U_a U_{-a}$ is trivial, we find $n = m(u) h$. As $m(u) \in \langle U_a, U_{-a} \rangle \cap N$, it follows that $h$ is also contained in this intersection. Hence $h$ stabilizes $x$ by the first case. It remains to show that $m(u)$ stabilizes $x$. As in the first case it follows that $m(u)$ centralizes all groups $U_b$ and $U_{-b}$ for $b \in Y$. Besides, $m(u)$ acts on $A$ as an affine reflection with invariant hyperplane $\{ z \in A: a(z) + s = 0\}$. If $a = \sum_{c \in [a]} n_c c$, we write $s = \sum_{c \in [a]} n_c s_c$ for suitable real numbers $s_c$. 
 Let $z_k$ be a sequence of points in $A$ satisfying $c(z_k) = c(x)$ for all $c \in Y$, $c(z_k) = -s_c$ for $c \in [a]$ and $c(z_k) = k$ for all $c \notin Y \cup [a]$. Then $a(z_k) +s = 0$, hence $m(u)$ stabilizes $z_k$.  Now let $\lambda$ be a $K$-weight such that $[\lambda_0(\Delta) - \lambda] \not\subset Y$. Then $[\lambda_0(\Delta) - \lambda]$ is not contained in $Y \cup [a]$, hence we find $(\lambda_0(\Delta) - \lambda)(z_k) \rightarrow \infty$. Therefore $z_k$ converges to $x$.
Since $N$ acts continuously on $\ov{A}$, we again deduce that $n$ stabilizes $x$. 
\hfill$\Box$

For all $\Omega \subset \ov{A}$ the group $N_\Omega$ normalizes $U_\Omega$ by Lemma 3.4. We put
\[P_\Omega = N_\Omega U_\Omega = U_\Omega N_\Omega.\]
Then $P_\Omega$ is a subgroup of $G$. 
If $\Omega = \{x\}$, we write $P_x = P_\Omega$.
\begin{cor} Let $\Omega$ be a non-empty subset of $\ov{A}$, and fix an order on $\Phi$ corresponding to some Weyl chamber.

i) $P_\Omega = U_\Omega^- U_\Omega^+ N_\Omega = N_\Omega U_\Omega^+ U_\Omega^-$

ii) $P_\Omega \cap U^\pm = U_\Omega^\pm$

iii) $P_\Omega \cap N = N_\Omega$.
\end{cor}
{\bf Proof: }i) By Theorem 3.6, we have $U_\Omega = U_\Omega^- U_\Omega ^+ (N \cap U_\Omega)$. Since $N \cap U_\Omega \subset N_\Omega$ by Corollary 3.7, our claim follows. 

ii) Obviously, $U_\Omega^-$ is contained in $P_\Omega \cap U^-$. By i), any $v$ in $P_\Omega \cap U^-$ can be written as $v = u^- u^+ n$ with $u^\pm \in U_\Omega ^\pm$
and $n \in N$. Then $n$ is contained in $U^+ U^- \cap N$, which is trivial by \cite{boti}, 5.15. Since $U^+ \cap U^- = \{1\}$, we find  $v = u^- \in U_\Omega^-$.

iii) Obviously, $N_\Omega$ is contained in $P_\Omega \cap N$. If $n \in P_\Omega \cap N$, we write it as $n = u^- u^+ m$ with $u^\pm \in U_\Omega^\pm$ and $m \in N_\Omega$.
Then $u^- u^+$ is contained in $N \cap U^- U^+ = \{1\}$, hence $n \in N_\Omega$.\hfill$\Box$

Now we can show a weak mixed Bruhat decomposition for the groups $P_x$, cf. \cite{we1}, Theorem 4.9, where this is proved in the case $\bg = {\bpgl}_n$. It is weak in the sense that one point has to lie in the apartment $A$. This makes the proof quite easy, since we can
deduce our statement from the Bruhat decomposition for two points in $A$.
This version suffices to define our compactification.
Possibly one can use the ideas in \cite{brti}, 7.3 to prove Iwasawa and Bruhat decompositions in general.

\begin{thm} Let $x \in \ov{A}$ and $y \in {A}$. Then $G = P_x N P_y$.
\end{thm}

{\bf Proof: }Let $A_F$ with  $F = \fyd$ be the stratum containing $x$. 
We choose a sequence of points $x_k \in A$ such that 
$a(x_k) = a(x)$ for all $a \in Y$ and such that $a(x_k) \rightarrow \infty$ for all $a \in \Delta\backslash Y$. Then $x_k$ converges towards $x$ in $\ov{A}$.

By the mixed Bruhat decomposition for points in $A$ (see \cite{brti}, 7.3.4) we have  $G = P_{x_k} N P_y$.

Consider $g\in G$. Now we apply corollary 3.8, where the order on $\Phi$ is defined by $\Delta$.  
For all $k$ we can write  $g = u_k^- u_k^+ n_k v_k$ for $u_k^\pm \in U_{x_k}^\pm$, $n_k \in N$ and $v_k \in P_y$. Since $P_y$ is a compact subgroup of $G$, we can pass to a subsequence and assume that $v_k$ converges towards some $v \in P_y$. 
For all $a \in \Phiredm$ the values
$-a(x_k)$ are bounded from below, hence $U_{a, x_k}$ is contained in some $U_{a,r}$ for $r \in \R$. Since  $U_{a,r}$ is a compact subgroup of $U_a$, there exists 
a compact subgroup of $U^-$ containing all 
$U_{x_k}^-$. By passing to a subsequence we can assume that $u_k^-$ converges
in the group $U^-$. Using Theorem 3.6 and  Lemma 3.5, we find that its limit $u^-$ lies in $U_x^-$. 

As $W = N/Z$ is a finite group, after passing once more to a subsequence we can write $n_k = z_k n$ for a sequence $z_k \in Z$ and some $n \in N$. Then $u_k^+ z_k$ is a convergent sequence in the minimal parabolic subgroup $U^+ Z$ associated to $\Delta$. Therefore its limit is equal to $u^+ z$ for some $u^+ \in U^+$ and $z \in Z$. Since the minimal parabolic $U^+ Z$ is the semidirect product of $U^+$ and $Z$, the sequence $u_k^+$ converges towards $u^+$, and the sequence $z_k$ converges towards $z$. Again, by Theorem 3.6 and Lemma 3.5, we find $u^+ \in U_x^+$. Passing to the limit in
$ g = u_k^- u_k^+ n_k v_k$ we get $g = u^- u^+ z n v \in U_x N P_y \subset P_x N P_y$.\hfill$\Box$

\section{Compactification of the building}

Now we use the groups $P_x$ investigated in the previous section to  define the compactified building.  We define a relation on $G \times \ov{A}$ by
\begin{eqnarray*}
(g,x) \sim (h,y) & \mbox{if and only if  there exists an element } n \in N \\
~ & \mbox{such that } nx= y \mbox{ and } g\inv h n \in P_x
\end{eqnarray*}
Since by Lemma 3.4  $n P_x n\inv = P_{nx}$ for all $n \in N$, this
is an equivalence relation. 

\begin{defi} 
The compactification of $X$ associated to the representation $\rho$  is defined as the quotient
\[\ov{X} = \ov{X}(\bg)_\rho = (G \times \ov{A}) / \sim.\]
\end{defi}

The space $\ov{X}$ carries a natural topology induced by the topology on $G \times \ov{A}$. 
The group $G$ acts continuously on $\ov{X}$ via multiplication in the first factor. 
Mapping $x \in \ov{A}$ to the class of $(1,x)$ in $\ov{X}$ defines a continuous map $\ov{A} \rightarrow \ov{X}$.
Since $N \cap P_x = N_x$ by Corollary 3.8, this is injective. 
Note that it follows directly from the definition, that for all points $x \in \ov{A}$ the group $P_x$ is the stabilizer of $x$. 

Now we want to show that $\ov{X}$ is compact. For this, we need the next result.

\begin{thm}
Let $(\gamma_k)_k$ be a sequence in $G$ converging to $\gamma \in G$ such that $\gamma_k \in N P_{x_k}$, where $x_k$ is a sequence in $\ov{A}$ converging to $x \in \ov{A}$. Then after passing to a subsequence of $(\gamma_k)_k$ we can write  $\gamma_k = n_k p_k$ with elements $n_k \in N$ and $p_k \in P_{x_k}$ such that $n_k$ converges to some $n \in N$ and $p_k$ converges to some $p \in P_x$. In particular, $\gamma$ lies in $N P_x$.
\end{thm}
{\bf Proof: }Recall that  $G$ acts continuously on $\ov{X}$ with stabilizer groups $P_x$.
Hence the claim is true if $\gamma_k \in P_{x_k}$. Thus, our task is to control the $N$-part of $\gamma_k$. 

Since there are only finitely many components $A_F$ in $\ov{A}$, we can pass to subsequences of $x_k$ and $\gamma_k$ and assume that all $x_k$ are contained in the same stratum $A_{F}$ for the face $F = \fyd$. Then the limit $x$ lies in $A_{F'}$ for some face $F'$ satisfying $F \subset \tilde{F}'$. 
After passing to a subsequence we can furthermore assume that all $x_k$ lie in the same Weyl chamber in $A_F$. Then there is an element $w \in W$ such that the induced map on $A_F$ maps the chamber in $A_F$ given by the base $Y$ to this new chamber. As in the proof of Proposition 2.4 we see that  $ \fyd  =  F_{w(Y)}^{w(\Delta)}$. Therefore we can replace $\Delta$ by $w(\Delta)$ and $Y$ by $w(Y)$ and assume that $b(x_k) \geq 0$ for all $b \in Y$. 

By corollary 3.8, we have $P_{x_k} = N_{x_k} U_{x_k}^+ U_{x_k}^-$, where we use the order on $\Phi$ given by $\Delta$. As $W = N/Z$ is a finite group, we can pass to subsequences and assume that there is an element $n  \in N$ such that $\gamma_k \in n Z U_{x_k}^+ U_{x_k}^-$ for all $k$. It suffices to prove our claim for $n\inv \gamma_k$, hence we can  assume that $\gamma_k \in Z U_{x_k}^+ U_{x_k}^-$.

Using Theorem 3.6, we write $\gamma_k = z_k w_k^+ u_k^- v_k^-$, with $z_k \in Z$, $w_k^+ \in U_{x_k}^+$, $u_k^- \in \prod_{a \in \Phiredm \backslash \langle Y \rangle} U_{a, x_k}$ and $v_k^- \in \prod_{a \in \Phiredm \cap \langle Y \rangle} U_{a, x_k}$.
For  $a \in \Phiredm \cap \langle Y \rangle$ we have $f_{x_k}(a) = -a(x_k) \geq 0 $, so that $U_{a,x_k}$ is contained in the compact subgroup $U_{a,0}$.
Hence  we can pass to a subsequence and assume that $v_k^-$ converges in $\prod_{a \in \Phiredm \cap \langle Y \rangle} U_{a}$. By Lemma 3.5 its limit lies in $U_x$. 

Put $Y^\perp = \{a \in \Delta: (a, b) = 0 \mbox{ for all } b \in Y\}$ and $\{\lambda_0(\Delta)\}^\perp = \{a \in \Delta: (\lambda_0(\Delta), a) = 0\}$, and look at the subset 
$X = Y^\perp \cap \{\lambda_0(\Delta)\}^\perp $ of $\Delta$.
If $a \in \Phiredm \backslash \langle Y \rangle$ such that $a \notin \langle  X \rangle$, then $[-a] \cup Y$ is an admissible  subset of $\Delta$. By the claim in the proof of Proposition 3.3 we find $U_{a, x_k} = 1$. Therefore $u_k^-$ lies in fact in $\prod_{a \in \Phiredm \cap  \langle X \rangle} U_{a, x_k}$. 

Hence we can replace $\gamma_k$ by $\gamma_k (v_k^-)\inv$ and assume that $\gamma_k = z_k w_k^+ u_k^-$ is contained in the subgroup of $G$ generated by $Z$ and all root groups $U_a$ for $a \in \langle X \rangle \cup \Phiredp$, i.e. in the standard parabolic ${\bp}_{X}$ associated to the subset $X$ of $\Delta$. By \cite{boti}, 4.2 and 5.12, ${\bp}_{X}$ is the semidirect product of  $Z(\bs_{X})$ and its unipotent radical $\prod_{a \in \Phiredp \backslash \langle X\rangle} \bu_a$, where $\bs_{X}$ is the identity component of the subgroup $(\bigcap_{a \in X} \mbox{ ker }a)$ of $\bs$. 
We write $w_k^+ = u_k^+ v_k^+$ with $u_k^+ \in \prod_{a \in \Phiredp \cap \langle X \rangle} U_{a,x_k}$ and $v_k^+ \in \prod_{a \in \Phiredp \backslash \langle X \rangle} U_{a, x_k}$. Then 
\[\gamma_k = z_k u_k^+ v_k^+ u_k^- = (z_k u_k^+ u_k^-) ((u_k^-)\inv v_k^+ u_k^-).
\]
The first factor converges in the Levi group $\bl = Z(\bs_{X})$,
the second factor converges in the unipotent radical $\prod_{a \in \Phiredp \backslash \langle X \rangle} \bu_a$. Since $(u_k^-)\inv v_k^+ u_k^-$ is contained in $P_{x_k}$, its limit lies in $P_{x}$. 

Therefore we can assume that the sequence $\gamma_k$ is of the form $\gamma_k = z_k u_k^+ u_k^-$, in particular it lies in $L = \bl(K)$. 
Note that $\bl$ is a reductive group over $K$ with maximal $K$-split torus $\bs$. 
Its root system $\Phi(\bs, \bl)$ can be identified with $\Phi \cap \langle X \rangle$. 

We look at the group $\bz$. By \cite{boti}, 4.28, there is a
 connected anisotropic subgroup $\bm$ of $\bz$ such that the multiplication map $\bm \times \bs \rightarrow \bz
$ is an isogeny. Its kernel $\bh$ can be identified with a finite subgroup of $\bs$. 
The obstruction for surjectivity of the map on $K$-rational points $M \times S \rightarrow Z$ is the first Galois cohomology of ${\bh}(\overline{K})$, which is a finite subgroup of a power of $\overline{K}^\times$. Since every finite subgroup of $\overline{K}^\times$ is cyclic of order prime to the characteristic of $K$, we deduce that also ${\bh}(\overline{K})$ has order prime to the characteristic of $K$ (if this is non-zero). Hence by \cite{se}, II,  5.2 its cohomology is finite.
Therefore we can pass to a subsequence and assume that all $z_k \in Z$ lie in the same coset with respect to the image of $M \times S$, i.e.
$z_k =  y m_k s_k$ with $y \in Z$,  $m_k \in M$ and  $s_k \in S$. 
Now $\bm$ is reductive and anisotropic, hence $M$ is compact by \cite{boti}, 9.4.
Passing to a subsequence, we can therefore assume that $m_k$ converges towards $m \in M$. Replacing $\gamma_k$ by $m_k\inv  y\inv \gamma_k$, we can assume that $z_k = s_k \in S$. 
 
Let $a \in  \Phired \cap \langle X \rangle$. Then for  every weight $\lambda$ with $[\lambda_0(\Delta) - \lambda] \subset Y$ and every $l \neq 0 $ the element $\lambda +la$ is not a weight.  Now we use our fixed $K$-rational representation $\rho: \bg \rightarrow {\bgl}(V)$. 
Put \[V_Y = \oplus_{\lambda: [\lambda_0(\Delta) - \lambda] \subset Y} V_\lambda,\]
where $V_\lambda$ is the weight space corresponding to $\lambda$.  
Then all root groups $\bu_a$ with $a \in \Phired \cap \langle X \rangle$ act trivially on $V_Y$ via $\rho$, which implies $\rho(\gamma_k)|_{V_Y} =  \rho(s_k)|_{V_Y}$. 

As $\gamma_k$ converges, the sequence $\rho(\gamma_k)|_{V_Y}$ converges in $GL(V_Y)$,
so that $\rho(s_k)|_{V_Y}$ converges in $GL(V_Y)$. 
Hence
for all weights $\lambda$ with $[\lambda_0(\Delta) - \lambda] \subset Y$ we find that
$\lambda(s_{k})$ converges. Therefore $a(s_{k})$ converges in $K^\times$ for all $a \in Y$. Let $[Y] \subset X^\ast(\bs)$ be the $\Z$-lattice generated by all $a \in Y$, and denote by $R$ the free part of the $\Z$-module $X^\ast(\bs) / [Y]$. Then there exists a $K$-split subtorus $\bs''$ of $\bs$ such that the corresponding map of character groups is the map $X^\ast(\bs) \rightarrow X^\ast(\bs'') = R$. Since $R$ is free, this map has a splitting. Therefore $\bs \simeq \bs' \times \bs''$ for some $K$-split subtorus $\bs'$ of $\bs$. 
Hence we can write $s_k = s_k' s_k''$ with $s_k' \in S'$ and $s_k'' \in S''$. By construction, the character group of $\bs'$ contains $[Y]$ as a subgroup of finite index. Hence for all characters $\chi$ in $X^\ast(\bs')$, the sequence $\chi(s_k')$ is bounded. Passing to a subsequence, we can assume that all $\chi(s_k')$ converge. Then $s_k'$ converges in $S'$. Therefore we can assume that $z_k = s_k'' \in S''$. Recall that $Z$ acts by translation with $\nu(z_k) = \nu(s_k'')$ on $A_F$. Since $ a(s_k'') = 1$ for all $a \in Y$, we find $\nu(z_k) = \nu(s_k'') \in \langle \fyd \rangle$, hence $z_k$ acts trivially on $A_F$. In particular, $z_k  \in Z_{x_k}$, i.e. $\gamma_k \in P_{x_k}$. 
Since $G$ acts continuously on $\ov{A}$, its limit $\gamma$ is then contained in $P_x$,
which finishes the proof.\hfill$\Box$

\begin{thm}
The compactification $\ov{X} = \ov{X}(\bg)_\rho$ of $X$ associated to the representation $\rho$ is a compact topological space.
\end{thm}
{\bf Proof: } Let $P_0$ denote the compact stabilizer group of the point $0$ in $A$. We restrict the equivalence relation $\sim$ to $P_0 \times \ov{A}$ and denote the quotient space as $ \ov{X}^\natural = (P_0 \times \ov{A}) / \sim$. Let us show that the relation $\sim$ is closed in $(P_0 \times \ov{A}) \times (P_0 \times \ov{A})$, which is a first countable space. Consider sequences  $g_k$, $h_k$ in $P_0$ and $x_k$, $y_k$  in $\ov{A}$ such that $(g_k, x_k) \sim (h_k, y_k)$ for all $k$. Assume that $g_k$ converges to $g$ and $h_k$ converges to $h$ in $P_0$, and also that $x_k$ converges to $x$ and $y_k$ converges to $y$ in $\ov{A}$.
We have to show that $(g,x) \sim (h,y)$. 

Since $(g_k, x_k) \sim (h_k,y_k)$, there exist elements $n_k \in N$ with $n_k x_k = y_k$ and $g_k\inv h_k n_k \in P_{x_k}$. Hence $ h_k\inv g_k$ lies in $n_k P_{x_k}$. By Theorem 4.2, we can write
$h_k\inv g_k = m_k f_k$ such that $m_k$ converges to $m \in N$ and $f_k$ converges to $f \in P_x$. Then $m_k\inv n_k$ lies in $N \cap P_{x_k} = N_{x_k}$. Hence we have $m_k x_k = n_k x_k = y_k$. Since $N$ acts continuously on $\ov{A}$ by Proposition 2.4, this implies $m x = y$. 
As $m_k f_k = h_k\inv g_k$ converges to $mf = h\inv g$, we find that $g\inv h m = f\inv \in P_x$.
Hence $(g,x) \sim (h,y)$, as claimed.

As $P_0 \times \ov{A}$ is a compact (i.e. Hausdorff and quasi-compact) space, it follows that the quotient space $\ov{X}^\natural$ after the closed equivalence relation $\sim$ is also compact.
Now we want to compare $\ov{X}$ and $\ov{X}^\natural$. The inclusion $P_0  \times \ov{A} \rightarrow G \times \ov{A}$ induces a continuous map $\alpha: \ov{X}^\natural \rightarrow \ov{X}$ on the quotient spaces. 
Besides, we define a map $\beta' : G \times \ov{A} \rightarrow \ov{X}^\natural$ as follows: Take $(g,x) \in G \times \ov{A}$
and write $g = pnq$ with $p \in P_0$, $n \in N$ and $q \in P_x$ using Theorem 3.9. Then let $\beta'(g,x)$ be the point in $\ov{X}^\natural$ induced by
$(p,nx) \in P_0 \times \ov{A}$. A straightforward calculation shows that this point does not depend on the choice of $p$, $n$ and $q$ in the decomposition of $g$, and that $\beta'$ factors over the equivalence relation $\sim$ on $G \times \ov{A}$. 

Let us now show that $\beta'$ is continuous. Consider a sequence $(g_k,x_k) \in G \times \ov{A}$ converging to $(g,x) \in G \times \ov{A}$. We write $g_k = p_k n_k 
q_k$ with $p_k \in P_0$, $n_k \in N$ and $q_k \in P_{x_k}$. After passing to a subsequence we can assume that $p_k$ converges towards $p$ in the compact group $P_0$.
Hence the sequence $n_k q_k \in N P_{x_k}$ converges in $G$. By Theorem 4.2, there is a sequence $m_k \in N$, converging towards $m \in N$ and a sequence $f_k \in P_{x_k}$ converging towards $f \in P_x$ such that $n_k q_k = m_k f_k$, i.e. $g_k = p_k m_k f_k$.  Then  $\beta'(g_k, x_k)$ is the point in $\ov{X}^\natural$ induced by $(p_k, m_k x_k)$. This sequence converges in $\ov{X}^\natural$ to the point induced by $(p, mx)$. Taking the limit of $g_k = p_k m_k f_k$, we find that $g = pmf$, so that $\beta'(g,x)$ is indeed the limit of $\beta'(g_k, x_k)$.  

Hence $\beta'$ induces a continuous map $\beta: \ov{X} \rightarrow \ov{X}^\natural$. By construction, $\alpha $ and $\beta$ are inverse to each other. Therefore
$\ov{X}$ is also compact.\hfill$\Box$

We will now show that the compactification $X(\bg)_\rho$  depends in fact only on the face of $X^\ast(\bs) \otimes \R$ containing the highest weight of $\rho$. In particular,
for every $\bg$ our construction yields a finite family of building compactifications.

First we go back to the definition of the compactified apartment and look at the faces $\fyd$ again. We want to  compare the faces $\fyd$ with the chamber faces given by the root system.
Let $\Delta$ be a base of the root system $\Phi$, and let $C(\Delta)$ denote the corresponding chamber in $A$, i.e.
$C(\Delta) = \{x \in A: a(x) > 0 \mbox{ for all } a \in \Delta\}$. For every $Z \subset \Delta$ we write 
\[E_Z = \{x \in A: a(x) = 0 \mbox{ for all }a \in Z \mbox{ and } a (x) > 0 \mbox{ for all } a \in \Delta \backslash Z\}\]
for the face corresponding to $Z$. Obviously, $E_\emptyset = C(\Delta)$. By $\tilde{C}(\Delta)$ we denote the closure of $C(\Delta)$ in $A$,
i.e. $\tilde{C}(\Delta) = \{ x \in A: a(x) \geq 0 \mbox{ for all } a \in \Delta\}$. 
\begin{prop}  Let $Y$ be an admissible subset of $\Delta$. 

i) Then
\[\fyd  \cap \tilde{C}(\Delta) = \bigcup ~_{Z: Y \subset Z \mbox{\footnotesize  max. adm. }} E_Z,\]
where the union on the right hand side runs over all subsets $Z$ of $\Delta$ such that $Y$ is the maximal subset of $Z$ which is admissible. In particular, $\fyd \cap \tilde{C}(\Delta)$ contains $E_Y$.

ii) Let $\Delta'$ be another base of $\Phi$. Then $\fyd \cap \tilde{C}(\Delta') \neq \emptyset$, if and only if there exists an admissible subset $Y'$ of $\Delta'$ with $\fyd = \fydp$.  
\end{prop}
{\bf Proof: } i) Let $x$ be a point in $E_Z$ for some $Z \subset \Delta$ containing $Y$ as a maximal admissible subset. Then $a(x) = 0$ for all $a \in Y$.
Since for any $K$-weight $\lambda$ the set $[\lambda_0(\Delta) - \lambda]$ is admissible, it is either contained in $Y$ or not contained in $Z$. In the second case,
$\lambda_0(\Delta) - \lambda$ is positive on $E_Z$, hence $E_Z \subset \fyd$.

Now we suppose $x \in \fyd \cap \tilde{C}(\Delta)$. Then $x$ lies in a face $E_Z$ for some $Z \subset \Delta$. Since $a(x) = 0$ for all $a \in Y$, we find
$Y \subset Z$. Assume that $Y'$ is an admissible subset of $Z$. Then there exists a $K$-weight $\lambda$ such that $[\lambda_0(\Delta) - \lambda] = Y'$. Since $Y' \subset Z$ and $x \in E_Z$, we have $(\lambda_0(\Delta) - \lambda)(x) = 0$. As $x \in \fyd$, this implies $Y' \subset Y$. Hence $Y$ is the maximal admissible subset of $Z$.

ii) Assume that there exists a point $x \in \fyd \cap \tilde{C}(\Delta')$. Then $x$ lies in some chamber face $E_{Z'} \subset \tilde{C}(\Delta')$. Let $Y'$ be the maximal admissible subset of $Z'$. Then
we have $x \in E_{Z'} \subset \fydp \cap \tilde{C}(\Delta')$ by i), which implies $x \in \fyd \cap \fydp$. By Lemma 2.1, it follows $\fyd = \fydp$. 
On the other hand, if $\fyd = \fydp$, then the intersection $\fyd \cap \tilde{C}(\Delta')$ is non-empty by i).\hfill$\Box$

Now we consider two faithful, geometrically irreducible rational representations $\rho$ and $\sigma$ of $\bg$ with highest weights $\lambda_0^{\rho}(\Delta)$ and $\lambda_0^\sigma(\Delta)$ (with respect to some base $\Delta$). Note that if there exists a homeomorphism $\varphi: \ov{X}(\bg)_\rho \rightarrow \ov{X}(\bg)_\sigma$ restricting to identity on $X(\bg)$, then it is uniquely determined by these properties. In this case we also say that $\ov{X}(\bg)_\rho$ and $\ov{X}(\bg)_\sigma$ are equal.

The root system $\Phi$ gives rise to finite collection of  hyperplanes $H_a = \{ z \in A^\ast: (a,z) = 0\}$ for $a \in \Phi$. We consider the associated decomposition of the dual space $A^\ast = X^\ast(\bs) \otimes \R$ into faces.

\begin{thm} The  compactifications $\ov{X}(\bg)_\rho$ and $\ov{X}(\bg)_\sigma$ are equal if and only if the highest weights $\lambda_0^\rho(\Delta)$ and $\lambda_0^\sigma(\Delta)$ are contained in the same face of $A^\ast$. 
\end{thm}
{\bf Proof: } Assume that there exists a homeomorphism $\varphi: \ov{X}(\bg)_\rho \rightarrow \ov{X}(\bg)_\sigma$  restricting to the identity map on $X(\bg)$. Since $\lambda_0^\rho(\Delta)$ and $\lambda_0^\sigma(\Delta)$ are dominant weights, we have $(\lambda_0^\rho(\Delta), a) \geq 0$ and $(\lambda_0^\sigma(\Delta),a) \geq 0$ for all $a \in \Delta$. Consider a root $a \in \Delta$ such that
$(\lambda_0^\rho(\Delta),a) > 0$. Then $Y = \{a\}$ is an admissible subset of $\Delta$ giving rise to a face $F = \fyd(\rho)$.  We choose a point $x \in A_F$ and a sequence of points $x_k \in A$ satisfying
$a(x_k) = a(x)$ and $b(x_k) \rightarrow \infty$ for all $b \in \Delta$ with $b \neq a$. Then $x_k$ converges to $x$ in $\ov{X}(\bg)_\rho$. Hence $x_k = \varphi(x_k) $ converges to $\varphi(x)$ in $\ov{X}(\bg)_\sigma$. We want to show that $(\lambda_0^\sigma(\Delta), a) > 0$. If $\lambda_0^\sigma(\Delta)$ is perpendicular to $a$, then $Y$ is not admissible with respect to $\sigma$. 
By definition of the sequence $x_k$, it follows that for all $K$-weights $\lambda \neq \lambda^\sigma_0(\Delta)$ the sequence $(\lambda_0^\sigma(\Delta) - \lambda)(x_k)$ converges to infinity, so that  $\varphi(x) \in A_{F'}$ for the face $F' = F_\emptyset^\Delta(\sigma)$. Now $A_{F'}$ consists of one point, which is by continuity the image of the one-point-set $A_{F_\emptyset^\Delta(\rho)}$. This contradicts injectivity. Hence $\lambda_0^\sigma(\Delta)$ cannot be perpendicular to $a$, i.e. $(\lambda_0^\sigma(\Delta),a) > 0$. Reversing the roles of $\rho$ and $\sigma$ we deduce that
$(\lambda_0^\rho(\Delta),a)  > 0$ if and only if $(\lambda_0^\sigma(\Delta),a) > 0$. Therefore $\lambda_0^\rho(\Delta)$ and $\lambda_0^\sigma(\Delta)$ lie in the same face.

In order to show the other direction, assume that $\lambda_0^\rho(\Delta)$ and $\lambda_0^\sigma(\Delta)$ lie in the same face. Then the graph structure on $\Delta \cup \{\larho\}$ coincides
with the graph structure of $\Delta \cup \{\lasigma\}$.
Hence a subset $Y$ of $\Delta$ is admissible with respect to $\lambda_0^\rho(\Delta)$ if and only if $Y$ is admissible with respect to $\lambda_0^\sigma(\Delta)$. 
Therefore we may call such a set $Y$ admissible, without specifying the representation. An admissible set $Y \subset \Delta$ gives rise to a face $\fyd(\rho)$ corresponding to $\rho$ and a face $\fyd(\sigma)$ corresponding to $\sigma$. We want to show that $\fyd(\rho) = \fyd(\sigma)$.

Let $Y'$ be an admissible subset of a base $\Delta'$. As a first step, we show the following claim:

{\bf Claim: }We have $\fyd(\rho) = \fydp(\rho)$ if and only if $\fyd(\sigma) = \fydp(\sigma)$.

By symmetry, it suffices to show one direction. Assume that $\fyd(\rho) = \fydp(\rho)$. By Lemma 2.1, this implies $\langle Y \rangle = \langle Y'\rangle$ and $\lambda_0^\rho(\Delta) - \lambda_0^\rho(\Delta') \in \langle Y \rangle$. There exists some $w$ in the Weyl group $W$ such that $w(\Delta) = \Delta'$, and hence $w(\lambda_0^\rho(\Delta)) = \lambda_0^\rho(\Delta')$.
By \cite{boti}, 12.17, we can write $w = w_1 w_2$ with $w_1, w_2 \in W$ such that $w_2(\lambda_0^\rho(\Delta)) = \lambda_0^\rho(\Delta)$ and such that 
$w_1 = s_{b_r} \cdots s_{b_1}$ is a product of reflections corresponding to roots $b_i \in [\lambda_0^\rho(\Delta)- \lambda_0^\rho(\Delta')]$. Since this
set is contained in $Y$, it follows that $w_1$ maps $\langle Y \rangle$ to itself.
Besides, since $w_2$ fixes $\lambda_0^\rho(\Delta)$, it fixes the whole face containing $\lambda_0^\rho(\Delta)$, in particular $\lambda_0^\sigma(\Delta)$. 
Therefore
\[
\lambda_0^\sigma(\Delta) - \lambda_0^\sigma(\Delta')  =\lasigma - w(\lasigma)
= \lasigma - w_1(\lasigma) ,\]
and a straightforward induction shows that this is contained in $\langle Y \rangle$. 
Since $\langle Y \rangle = \langle Y' \rangle$, this implies 
$\fyd(\sigma) = \fydp(\sigma)$ by Lemma 2.1, which proves the claim.

Now we show $\fyd(\rho)= \fyd(\sigma)$ for all admissible subsets $Y$ of a base $\Delta$. It suffices to show that for all bases $\Delta'$ of $\Phi$ we have
$\fyd(\rho) \cap \tilde{C}(\Delta') = \fyd(\sigma) \cap \tilde{C}(\Delta')$. 
It follows from Proposition  4.4 combined with the claim, that  $\fyd(\rho) \cap \tilde{C}(\Delta')$ is non-empty if and only if  $\fyd(\sigma) \cap \tilde{C}(\Delta') $ is non-empty.
If both $\fyd(\rho) \cap \tilde{C}(\Delta')$ and $\fyd(\sigma) \cap \tilde{C}(\Delta')$ are non-empty, then by Proposition 4.4, $\fyd(\rho) = \fydp(\rho)$ for some admissible subset $Y'$ of $\Delta'$. We deduce from the claim that
$\fyd(\sigma) = \fydp (\sigma)$. Therefore by Proposition 4.4, the sets $\fyd(\rho) \cap \tilde{C}(\Delta')$ and $\fyd(\sigma) \cap \tilde{C}(\Delta')$ are both the union of all chamber faces $E_{Z'}$, where $Z'$ is a subset of $\Delta'$ containing $Y'$ as a maximal admissible subset. Hence they are equal.

Therefore $\rho$ and $\sigma$ give rise to the same decomposition of $A$ into faces $\fyd$. Hence the compactified apartments for $\rho$ and $\sigma$ coincide. Besides, the groups $P_x$ for $x \in \ov{A}$ coincide, so that $\ov{X}(\bg)_\rho$ is indeed equal to $\ov{X}(\bg)_\sigma$.\hfill$\Box$

\section{Comparison with other compactifications}
In this section we show that the compactifications defined in \cite{la}, \cite{we1} and \cite{we2}
occur among the compactifications studied in the present paper.
However, the Borel-Serre compactification defined in \cite{bose} is of a completely different nature and does not occur in our family. 

First  we consider the polyhedral compactification $\ov{X}(\bg)^P$ sketched in \cite{ge} and constructed in \cite{la}.
Let $\bg$ be a connected, semisimple group. If $\rho$ is a geometrically irreducible $K$-rational representation 
with heighest weight lying in the (open) chamber of $A^\ast$ associated to the base $\Delta$, then
any subset $Y$ of $\Delta$ is admissible. Hence the decomposition of $A$ into the faces $\fyd$ coincides with the decomposition of $A$ into the chamber faces given by the hyperplanes
$\{x \in A: a(x) = 0\}$ for all $a \in \Phi$. Therefore our compactification $\ov{A}$ is naturally homeomorphic in an $N$-equivariant way to the compactified apartment $\ov{A}^P$ in the polyhedral compactification of $X$, see \cite{la}, chapter I. Since our groups $P_x$ are defined in a similar way as in \cite{la}, 12.4, this homeomorphism can be continued to a natural $G$-equivariant homeomorphism $\ov{X}(\bg)_\rho \iso \ov{X}(\bg)^P$.  
 
In \cite{we1} and \cite{we2} we defined two compactifications of the building $X({\bpgl} (V))$. In order to fit these into the frame of the present paper, we replace ${\bpgl} (V)$ by ${\bsl} (V)$. The underlying topological space of the building is the same in both cases, and the ${\bpgl}(V)$ action on $X({\bpgl} (V)) = X({\bsl}(V))$ is induced from the ${\bsl}(V)$-action.

Let $i: {\bsl}(V) \rightarrow {\bgl}(V)$ be the natural embedding. Then it is easily seen that the compactified apartment $\ov{A}$ with respect to $i$ (together with the $N$-action) coincides with the compactified apartment in \cite{we1} (together with the $N$-action). Therefore the groups $P_x$ defined here coincide with the corresponding groups in \cite{we1}. Hence  it follows that the compactification $\ov{X}({\bsl}(V))_i$ coincides with the compactification in \cite{we1}.

The same reasoning shows that the compactification defined in \cite{we2} using seminorms on $V$ coincides with the compactification $\ov{X}({\bsl}(V))_{i^\ast}$ for the dual representation $i^\ast$.

\section{Further questions}
We conclude this paper with some open questions which will be the subject of further research.

The finite familiy of compactifications $\ov{X}(\bg)_\rho$ is in some respects analogous to the different Satake compactifications of a symmetric space, see \cite{sa}.  Let ${\cal S} = G/K$ be a Riemann symmetric space of non-compact type, and fix a fatihful representation of $G$ in $PSL(n,\C)$ satisfying the conditions in \cite{sa}, 2.2. The symmetric space associated to $PSL(n,\C)$ has a natural compactification $\cal Z$, which is the projective space associated to the space of all hermitian $n \times n$-matrices. Via the representation, $\cal S$ can be embedded in $\cal Z$.  Taking the closure of $\cal S$ in $\cal Z$ yields a compactification. The structure of its boundary is investigated in \cite{sa}, section 4.

 Is it possible to describe our building compactifications $\ov{X}(\bg)_\rho$ in a similar way? The space $\cal Z$ looks very much like the Archimedean analogue of the seminorm compactification for $X({\bpgl})$ defined in \cite{we2}.  Here one has to take into account that functoriality of Bruhat-Tits buildings is a delicate question, see e.g. \cite{la2}.

In the family of Satake compactifications, there is a maximal one, which dominates the others. This maximal Satake compactification is analogous to the polyhedral compactification for buildings.
Presumably also in the building case, the polyhedral compactification is maximal in the sense that it dominates all the other family members. 

Besides, it seems likely that the boundary of $X(\bg)_\rho$ can be identified with the union
of buildings corresponding to the parabolic subgroups of $\bg$ of a certain type. This result should be analogous to the description of the boundary components in the Satake compactification given in \cite{sa}, Theorem 1. 

The compactification of $X({\bpgl}(V))$ defined in \cite{we2} can be identified with 
a closed subset of a certain Berkovich topological space, see \cite{we2}, Proposition 6.1. In view of the connection between Bruhat-Tits buildings and Berkovich spaces (see \cite{ber}, chapter 5), it seems 
possible to identify all our compactifications $\ov{X}(\bg)_\rho$ with subsets of suitable Berkovich spaces. 

The recent paper \cite{guire} by Guivarc'h and R\'emy raises another question: Is there a group-theoretic compactification of the vertex set in the building which can be identified with the closure of the vertex set in $\ov{X}(\bg)_\rho$?

\end{document}